\documentclass[journal]{IEEEtran}

\usepackage{cite}
\usepackage{url}
\usepackage[caption=false]{subfig}
\usepackage{float}
\usepackage[hidelinks]{hyperref}
\usepackage{amsmath,amssymb,amsfonts}
\allowdisplaybreaks[4]
\usepackage{graphicx}
\usepackage{enumerate}
\usepackage{array}
\usepackage{multirow}
\usepackage{color}
\usepackage[ruled,linesnumbered]{algorithm2e}

\usepackage{tikz}

\newtheorem{theorem}{Theorem}
\newtheorem{lemma}{Lemma}

\begin{document}

\title{Predict-and-Optimize Robust Unit Commitment with Statistical Guarantees via Weight Combination}

\author{Rui~Xie,~\IEEEmembership{Member,~IEEE,}
        Yue~Chen,~\IEEEmembership{Senior Member,~IEEE,~}
        and Pierre Pinson,~\IEEEmembership{Fellow,~IEEE}
\thanks{R. Xie and Y. Chen are with the Department of Mechanical and Automation Engineering, The Chinese University of Hong Kong, HKSAR, China. (email: ruixie@cuhk.edu.hk; yuechen@mae.cuhk.edu.hk) (Corresponding author: Yue Chen)}
\thanks{P. Pinson is with the Dyson School of Design Engineering, Imperial College London, UK. (email: p.pinson@imperial.ac.uk)}
}

\markboth{Journal of \LaTeX\ Class Files,~Vol.~14, No.~8, August~2024}%
{Shell \MakeLowercase{\textit{et al.}}: Bare Demo of IEEEtran.cls for IEEE Journals}

\maketitle

\begin{abstract}
The growing uncertainty from renewable power and electricity demand brings significant challenges to unit commitment (UC). While various advanced forecasting and optimization methods have been developed to better predict and address this uncertainty, most previous studies treat forecasting and optimization as separate tasks. This separation can lead to suboptimal results due to misalignment between the objectives of the two tasks. To overcome this challenge, we propose a robust UC framework that integrates forecasting and optimization processes while ensuring statistical guarantees. In the forecasting stage, we combine multiple predictions derived from diverse data sources and methodologies for an improved prediction, aiming to optimize the UC performance. In the optimization stage, the combined prediction is used to construct an uncertainty set with statistical guarantees, based on which the robust UC model is formulated. 
The optimal robust UC solution provides feedback to refine the weight used for combining multiple predictions. 
To solve the proposed integrated forecasting-optimization framework efficiently and effectively, we develop a neural network-based surrogate model for acceleration and introduce a reshaping method for the uncertainty set based on the optimization result to reduce conservativeness. Case studies on modified IEEE 30-bus and 118-bus systems demonstrate the advantages of the proposed approach.
\end{abstract}



\section*{Nomenclature}

\subsection{Abbreviation}

\begin{IEEEdescription}[\IEEEusemathlabelsep\IEEEsetlabelwidth{ssssssssss}]
\item[C\&CG] Column-and-constraint generation.
\item[DRO] Distributionally robust optimization.
\item[i.i.d.] Independent and identically distributed.
\item[MILP] Mixed-integer linear programming.
\item[MLP] Multilayer perceptron.
\item[PSO] Particle swarm optimization.
\item[RO] Robust optimization.
\item[SP] Stochastic programming.
\item[UC] Unit commitment.
\end{IEEEdescription}

\subsection{Indices and Sets}

\begin{IEEEdescription}[\IEEEusemathlabelsep\IEEEsetlabelwidth{ssssssssss}]
\item[$i \in \mathcal{I}$] Set of buses.
\item[$m \in \mathcal{M}$] Set of prediction methods.
\item[$x \in \mathcal{X}$] Feasible set of the pre-dispatch variable.
\item[$u \in \mathcal{U}$] Uncertainty set.
\item[$y \in \mathcal{Y}(x, u)$] Feasible set of the re-dispatch variable.
\item[$g \in \mathcal{G}$] Set of generators.
\item[$t \in \mathcal{T}$] Set of periods in unit commitment.
\item[$l \in \mathcal{L}$] Set of transmission lines.
\end{IEEEdescription}
       
\subsection{Parameters}

\begin{IEEEdescription}[\IEEEusemathlabelsep\IEEEsetlabelwidth{ssssssssss}]
\item[$T$] The number of periods in a day.
\item[$\varepsilon, \delta$] Probability tolerance parameters.
\item[$\hat{u}_{i t}$] Prediction of load $i$ in period $t$.
\item[$o_g^+, o_g^-$] Startup/shutdown cost of generator $g$.
\item[$\rho_g$] Generation cost coefficient of generator $g$.
\item[$\gamma_g^+, \gamma_g^-$] Upward/downward reserve cost coefficient of generator $g$.
\item[$S_l$] Capacity of line $l$.
\item[$\pi_{g l}, \pi_{i l}$] Power transfer distribution factor from generator $g$/bus $i$ to line $l$.
\item[$T_g^+, T_g^-$] Minimum up/down time of generator $g$.
\item[$R_g^+, R_g^-$] Maximum upward/downward reserve of generator $g$.
\item[$\underline{P}_g, \overline{P}_g$] Minimum/maximum output of generator $g$.
\item[$K_g^+, K_g^-$] Maximum upward/downward ramp of generator $g$ in a period.
\item[$K_g^U, K_g^D$] Maximum output increase/decrease when generator $g$ startups/shutdowns in a period.
\item[$\rho_g^+, \rho_g^-$] Upward/downward output adjustment cost coefficient of generator $g$.
\item[$\underline{U}, \overline{U}$] Lower/upper bound of uncertain load.
\end{IEEEdescription}

\subsection{Variables}

\begin{IEEEdescription}[\IEEEusemathlabelsep\IEEEsetlabelwidth{ssssssssss}]
\item[$u_{i t}$] Uncertain load at bus $i$ in period $t$.
\item[$w$] Weight vector of predictions.
\item[$\theta_{g t}, \theta_{g t}^+, \theta_{g t}^-$] Indicator variable for the on/startup/shutdown state of generator $g$ in period $t$.
\item[$p_{g t}$] Day-ahead scheduled output of generator $g$ in period $t$. 
\item[$r_{g t}^+, r_{g t}^-$] Upward/downward reserve of generator $g$ in period $t$.
\item[$p_{g t}^+, p_{g t}^-$] Upward/downward output adjustment of generator $g$ in period $t$.
\end{IEEEdescription}

\section{Introduction}

\IEEEPARstart{T}{he} ongoing transition towards greener power systems has significantly increased the deployment of renewable power generators, leading to higher volatility in power sources. This uncertainty, together with the randomness of electric loads, presents great challenges for power system operations. Unit commitment (UC) is one of such power system operation problems that require particular attention.

If given the probability distribution of uncertainty, stochastic programming (SP) can be applied to UC problems to determine the optimal strategy \cite{paturet2020stochastic}. However, an accurate probability distribution is difficult to obtain in reality and the inaccuracy of the distribution may jeopardize the feasibility and optimality of SP solutions. Robust optimization (RO) handles this problem by focusing on the worst-case scenario in a pre-defined uncertainty set. Fruitful research has been done in this area: A multi-stage robust UC model was proposed and solved by robust dual dynamic programming in \cite{xiong2022multi}. Polyhedra were used to replace rectangular uncertainty sets to model the uncertainty correlations in \cite{zhou2023partial}. However, RO can be overly conservative when the uncertainty set contains all the possible uncertainty realizations, because extreme scenarios rarely happen but can cause a steep cost increase. The method in \cite{ju2023two} estimated quantiles from historical data to establish rectangular uncertainty sets, in which extreme scenarios were excluded. In \cite{wang2024two}, uncertainty budgets were used to restrict the uncertainty set. A decision-dependent uncertainty set was constructed in \cite{haghighat2024robust} to consider the impact of pricing strategies on demand. However, estimated quantiles and uncertainty budgets may still be inaccurate, and the resulting uncertainty set lacks a statistical guarantee. This makes the obtained robust UC strategy less reliable and trustworthy. 


Distributionally robust optimization (DRO) is another way to deal with uncertainty, which considers the worst-case probability distribution in a pre-defined ambiguity set. Two-stage distributionally robust UC methods were developed in \cite{wang2022distributionally} and \cite{zheng2023day} using Wasserstein metric ambiguity sets, while unimodality skewness of wind power was utilized in \cite{zhou2023distributionally}. The copula theory was combined with DRO in \cite{liu2024modeling} to better capture the dependence between uncertainty. 
While data-driven DRO can offer statistical guarantees by constructing the ambiguity set based on the asymptotic properties of the empirical probability distribution \cite{wang2022distributionally}, the theoretical number of required data points typically depends on the dimension of uncertainty \cite{gao2023finite}, which is often much larger than what is practically available under multi-dimensional uncertainty, such as in the UC problem. Consequently, the size of the ambiguity set is typically adjusted through a trial-and-error process \cite{xie2022sizing}, without providing a satisfactory statistical guarantee.

Recently, a data-driven uncertainty set construction framework was proposed in \cite{hong2021learning} with a dimension-free statistical guarantee. Even for multi-dimensional uncertainty, it can provide a satisfactory statistical guarantee based on a moderate amount of data, enabling the decision-maker to effectively control the conservative degree and easily strike a balance between optimality and robustness. 
Therefore, this method is more suitable than DRO for the UC problem, where the dimension of uncertainty is proportional to the product of the number of periods and the number of uncertain net loads.
This method was applied to economic dispatch \cite{lu2024sample}, scheduling of thermostatically controlled loads \cite{jiang2023robust}, and UC \cite{liang2024joint}. However, the static RO models in \cite{lu2024sample,jiang2023robust,liang2024joint} cannot account for the re-dispatch stage, resulting in overly conservative day-ahead strategies. Therefore, a two-stage adjustable RO method with statistical guarantees is needed for the robust UC problem.

In addition, the quality of the uncertainty set relies heavily on the quality of the uncertainty forecasts. The predictions of different forecasting methods can be combined by a weight parameter to achieve a better performance than the individual methods \cite{wang2018ensemble}. Following this idea, an ensemble deep learning method was proposed in \cite{cui2022ensemble} for probabilistic wind power forecasting, where the weight parameter was determined according to a quantile loss index. The extreme prediction risk was modeled by the conditional value-at-risk in \cite{wang2023risk} to optimize the ensemble weight of renewable energy forecasts, and then the risk of renewable energy bidding strategy was evaluated. To enhance the performance in a changing environment, a deep deterministic policy gradient-based method was proposed in \cite{li2023adaptive} to adjust the combination weight adaptively. The ensemble forecasting framework was integrated with flexible error compensation in \cite{su2024towards}, and the weight was optimized to minimize the worst-case forecast error. The above studies chose the weights to optimize the accuracy of the forecasts, without considering the impact of the forecasts on the subsequent decision-making.

Conventionally, forecasting and optimization are performed separately, with forecasting focused on maximizing prediction accuracy and optimization aimed at minimizing costs. However, since the objectives of these two processes are distinct and may even conflict, conducting them separately can lead to suboptimal outcomes. 
To overcome this shortcoming, a predict-and-optimize framework was proposed in \cite{elmachtoub2022smart}, which integrates the forecasting and optimization processes and aims at improving the performance of the final strategy. 

The predict-and-optimize framework has been applied in various power system applications. For instance, in \cite{dias2025application}, the convergence of estimators derived from closed-loop prediction and optimization was examined in the context of linear programming. Solution methods utilizing MILP and heuristic algorithms were introduced and applied to reserve sizing in power systems. In \cite{chen2022feature}, uncertainty prediction and network-constrained UC were integrated into a closed-loop predict-and-optimize framework, with the prediction model trained according to the induced UC cost. The predict-and-optimize framework was also employed in \cite{sang2022safety} for voltage regulation in distribution systems. In \cite{wu2024novel}, stochastic UC was incorporated into the predict-and-optimize framework, where parameters in wind power forecasting models were adjusted based on risk costs using an iterative solution algorithm. Ref. \cite{shen2024predict} combined the predict-and-optimize framework with distributionally robust UC, adjusting the inertia prediction model according to UC cost feedback, while addressing price, demand, and renewable uncertainty through a scenario-based DRO technique. An end-to-end approach for multi-energy load data valuation was proposed in \cite{zhou2024load}, where the model was trained through forward and backward propagation of the forecasting model and the optimization problem. Ref. \cite{zhuang2025weighted} introduced a novel predict-and-optimize framework that considers the relative importance of uncertainty at different nodes, applying it to the optimal power flow problem in distribution networks.

However, there are still research gaps in the aforementioned predict-and-optimize framework: First, it directly adjusts the parameters of the forecasting model based on the induced cost from the optimized strategy, which may lead to biased predictions with weak explainability and transparency. Second, the presence of numerous adjustable parameters in the forecasting model results in a significant computational burden during parameter optimization within the predict-and-optimize framework. Third, to the best of the authors' knowledge, the predict-and-optimize framework has not been combined with a statistical guarantee, which could theoretically ensure the out-of-sample performance of the optimized strategy.

To bridge the aforementioned research gaps, this paper proposes a novel UC framework that integrates the forecasting and optimization processes while ensuring statistical guarantees. The main contribution is two-fold:


1) An integrated forecasting and optimization framework is proposed to predict in a way that optimizes the performance of the robust strategy. 
Specifically, predictions from various sources are combined using weight parameters to generate an improved forecast for RO; and in turn, the weights are optimized based on the evaluated performance of the resulting strategy. 
To accelerate the weight optimization process, a neural network is trained as a surrogate model and equivalently transformed into mixed-integer linear constraints. By solving a mixed-integer linear programming (MILP) problem, the optimal weights can then be obtained. Notably, this integrated framework and solution methodology have not been previously reported in the literature.


2) A data-driven two-stage robust UC model with statistical guarantees is proposed. It employs the prediction and historical data to construct an ellipsoidal uncertainty set, which is then enhanced by reconstructing a polyhedral uncertainty set based on the identified feasible solution and the UC parameters. This method provides dimension-free statistical guarantees to ensure robustness at a specified confidence level, which extends the method in \cite{hong2021learning} for static RO to two-stage RO.

The rest of this paper is organized as follows: The integrated forecasting and optimization framework is proposed in Section~\ref{sec:framework}. A two-stage robust UC method with statistical guarantees is developed in Section~\ref{sec:RUC}. Case studies are presented in Section~\ref{sec:case} and conclusions are drawn in Section~\ref{sec:conclusion}.

\section{Integrated Forecasting and Optimization Framework}
\label{sec:framework}

In this section, we first propose an integrated forecasting and optimization framework for a two-stage decision-making problem under uncertainty. In fact, the proposed framework is generally applicable and not restricted to the UC problem. Then in Section \ref{sec:RUC}, we introduce the detailed UC model and explain how it fits within the proposed framework. 


\subsection{Forecasting}
\label{sec:framework-forecasting}

Consider the forecasting task for a multivariate stochastic process $\mathbf{u} := (\mathbf{u}_1, \mathbf{u}_2, \dots, \mathbf{u}_t, \dots)$, where $\mathbf{u}_t = (u_{i t}; i \in \mathcal{I})$ is the vector of uncertain variables in period $t$, and $\mathcal{I}$ is the index set. A forecast is then a collection of future values over the horizon of interest, denoted by $\hat{\mathbf{U}} := (\hat{\mathbf{u}}_1, \hat{\mathbf{u}}_2, \dots, \hat{\mathbf{u}}_T)$. The hat notation indicates these are forecasts, and hence estimates of such future values. Additionally, the corresponding ground truth is denoted by $\mathbf{U}$.

Multiple forecasting methods can be used simultaneously. The prediction by method $m \in \mathcal{M}$ is denoted by $\hat{\mathbf{U}}^{(m)}$, and $\mathcal{M}$ is the set of methods. The final prediction is a convex combination of these predictions, that is, $\hat{\mathbf{U}}(w) = \sum_{m \in \mathcal{M}} w^{(m)} \hat{\mathbf{U}}^{(m)}$, where $w = (w^{(m)}; m \in \mathcal{M})$ is the weight vector. It satisfies $w \in \mathcal{W} := \{ w ~|~ w \geq 0, \sum_{m \in \mathcal{M}} w^{(m)} = 1\}$. Using the ground truth $\mathbf{U}$, the forecast error is calculated as $\mathbf{E}(w) = \mathbf{U} - \hat{\mathbf{U}}(w)$, which also depends on the weight $w$. We will denote the collection of $N$ such forecast error groups as $\mathbf{E}^{(1: N)}(w)$. Forecasting methods are not the focus of this paper, and we apply multiple forecasting methods in case studies based on the related literature.

\subsection{Two-Stage Optimization Under Uncertainty}
\label{sec:framework-optimization}

The following is a two-stage optimization formulation under uncertainty:
\begin{subequations}
\label{eq:UC-original}
\begin{align}
    \label{eq:UC-original-1}
    O := \min\nolimits_{x \in \mathcal{X}, \eta}~ & f(x) + \eta \\
    \label{eq:UC-original-2}
    \mbox{s.t.}~ & \Pr \left[ \min\nolimits_{y \in \mathcal{Y}(x, u)} h(y) \leq \eta \right] \geq 1 - \varepsilon.
\end{align}
\end{subequations}
In \eqref{eq:UC-original}, the first-stage variables are collected in $x$ and $\mathcal{X}$ is its feasible region. The first-stage cost function is denoted by $f(x)$. The uncertainty is denoted by $u$, the second-stage variables are collected in $y$, $\mathcal{Y}(x, u)$ is the feasible region of $y$ depending on $x$ and $u$, and $h(y)$ is the second-stage cost function, which is minimized after the realization of $u$ is observed. 
$\Pr[\cdot]$ denotes the probability. 
The second-stage cost value is denoted by $\eta$, and according to \eqref{eq:UC-original-2}, the probability of the second-stage cost being no larger than $\eta$ is at least $1 - \varepsilon$, where $\varepsilon$ is a specified probability threshold. 
Therefore, $\eta$ represents a $(1 - \varepsilon)$-quantile of the second-stage cost, and problem \eqref{eq:UC-original} minimizes the total cost.

Although problem \eqref{eq:UC-original} is a chance-constrained SP problem that has been studied and applied extensively in the literature, there is a key difficulty in its solution, i.e., the accurate probability distribution of uncertainty $u$ is usually unknown. In the case considered, we need to extract the distribution information from the historical data. On the one hand, the number of historical data is limited and the error of the empirical distribution is inevitable. On the other hand, we want to ensure the robustness, or more specifically, guarantee the chance constraint \eqref{eq:UC-original-2}. For the sake of robustness, instead of directly using the empirical distribution in problem \eqref{eq:UC-original}, we construct an uncertainty set $\mathcal{U}$ subject to $\Pr[u \in \mathcal{U}] \geq 1 - \varepsilon$, and consider the following two-stage RO problem:
\begin{align}
    \label{eq:RUC-original}
    O_{\mathcal{U}} := \min\nolimits_{x \in \mathcal{X}} \left\{ f(x) + \max\nolimits_{u \in \mathcal{U}} \min\nolimits_{y \in \mathcal{Y}(x, u)} h(y) \right\}.
\end{align}
The effectiveness of problem \eqref{eq:RUC-original} is revealed in Lemma~\ref{lemma:RUC}.

\begin{lemma}
    \label{lemma:RUC}
    Suppose $(x^*, \eta^*)$ and $x_\mathcal{U}^*$ are optimal solutions to problems \eqref{eq:UC-original} and \eqref{eq:RUC-original} with optimal values $O$ and $O_\mathcal{U}$, respectively. For any solution $x \in \mathcal{X}$, let
    \begin{align}
        O_x := f(x) + \min \left\{ \eta ~\middle|~ \Pr \left[ \min_{y \in \mathcal{Y}(x, u)} h(y) \leq \eta \right] \geq 1 - \varepsilon \right\}. \nonumber
    \end{align}
    Then $O = O_{x^*} \leq O_{x_\mathcal{U}^*} \leq O_\mathcal{U}$ whenever $\mathcal{U}$ satisfies $\Pr[u \in \mathcal{U}] \geq 1 - \varepsilon$. 
    Moreover, $O = O_{\mathcal{U}^*}$ and $\Pr[u \in \mathcal{U}^*] \geq 1 - \varepsilon$ for
    \begin{align}
        \label{eq:U*}
        \mathcal{U}^* := \left\{ u ~\middle|~ \min\nolimits_{y \in \mathcal{Y}(x^*, u)} h(y) \leq \eta^* \right\}. 
    \end{align}
\end{lemma}

Lemma~\ref{lemma:RUC} shows that given $\Pr[u \in \mathcal{U}] \geq 1 - \varepsilon$, the RO problem \eqref{eq:RUC-original} is a conservative approximation of problem \eqref{eq:UC-original}, and it becomes exact when $\mathcal{U} = \mathcal{U}^*$. Moreover, the performance of the obtained solution $x_\mathcal{U}^*$ is represented by $O_{x_\mathcal{U}^*}$, which is bounded from above by the optimal value $O_\mathcal{U}$ of problem \eqref{eq:RUC-original}. The proof of Lemma~\ref{lemma:RUC} is in Appendix~\ref{appendix:lemma-RUC}. 

The choice of the uncertainty set $\mathcal{U}$ is the key to achieving good performance in problem \eqref{eq:RUC-original}, for which we will leverage the information we have: Given a fixed $w$, the uncertainty set $\mathcal{U}(w)$ will be constructed using the prediction $\hat{\mathbf{U}}(w)$ and its accuracy estimation deduced from the historical data of forecast error $\mathbf{E}^{(1: N)}(w)$ to approach $\Pr [u \in \mathcal{U}(w)] \geq 1 - \varepsilon$. Subsequently, RO in \eqref{eq:RUC-original} is solved to find the optimal first-stage strategy $x_\mathcal{U}^*(w)$. The details of uncertainty set construction and solution method are in Section~\ref{sec:RUC}.

\subsection{Performance Evaluation and Integrated Framework}
\label{sec:framework-integration}

The proposed framework aims to solve the chance-constrained problem \eqref{eq:UC-original} to minimize costs while ensuring robustness. Problem \eqref{eq:UC-original} minimizes the summation of the first-stage cost and the $(1 - \varepsilon)$-quantile of the second-stage cost. Thus, quantile estimation is required to evaluate a solution. In daily operational applications like UC, problem \eqref{eq:UC-original} is solved daily, but only one actual observed data point is available per day, which is insufficient for quantile estimation. To address this, we approximate the $(1 - \varepsilon)$-quantile in cost evaluation using the $\lceil (1 - \varepsilon)N' \rceil$-th smallest cost from the $N'$ uncertainty data constructed based on the prediction and the historical data of forecast errors. 

Specifically, the performance of the optimized strategy $x_\mathcal{U}^*(w)$ is evaluated using the second training dataset, which is independent of the first training dataset used to construct $\mathcal{U}$.
The evaluated cost will be utilized to adjust $w$ in $\mathcal{W}$. Denote the historical data of forecast errors in the second training dataset by $\mathbf{E}^{(N + 1: N + N')}(w)$. The uncertainty realization data $\check{\mathbf{U}}^{(1: N')}$ for cost evaluation is constructed by $\check{\mathbf{U}}^{(n)} = \hat{\mathbf{U}}(w) + \mathbf{E}^{(N + n)}(w)$, $n = 1, 2, \dots, N'$. 
For each $\check{\mathbf{U}}^{(n)}$, the total cost is 
\begin{align}
    I^{(n)}(\hat{\mathbf{U}}, w) := f(x_\mathcal{U}^*(w)) + \min\nolimits_{y \in \mathcal{Y}(x_\mathcal{U}^*(w), \check{\mathbf{U}}^{(n)})} h(y). \nonumber
\end{align}
Sort $I^{(1: N')}(\hat{\mathbf{U}}, w)$ in an ascending order. The evaluated cost is then determined as the $\lceil (1 - \varepsilon)N' \rceil$-th cost (where $\lceil \cdot \rceil$ means rounding up to the nearest integer), denoted by $I(\hat{\mathbf{U}}, w)$.

\begin{figure}[!t]
   \centering
   \includegraphics[width=0.48\textwidth]{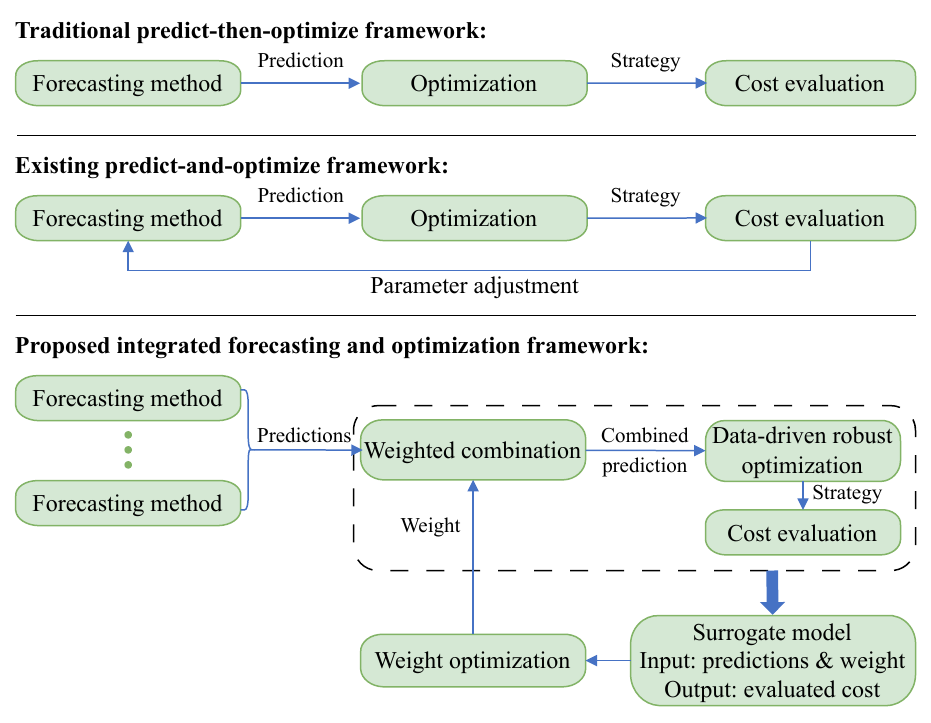}
   \caption{Forecasting and optimization frameworks.}
   \label{fig:framework}
\end{figure}

The proposed integrated forecasting and optimization framework is illustrated and compared with other approaches in Fig.~\ref{fig:framework}. The predictions from various forecasting methods are combined using the weight $w$. The combined prediction is then utilized in data-driven RO, where statistical guarantees are provided through the construction of uncertainty set. We propose to tune the weight $w$ in $\mathcal{W}$ to minimize the evaluated cost of the strategy. For applications like UC, optimizing $w$ must be efficient. Therefore, we develop a surrogate model in the next subsection that captures the relationship between the inputs (predictions and weight) and the output (evaluated cost). This surrogate model is trained offline. For daily optimization, the surrogate model is utilized to establish an optimization problem that optimizes the weight and minimizes the evaluated cost. This approach accelerates the weight optimization process. The optimized weight is then applied in the data-driven RO to derive the final strategy.

\subsection{Surrogate Model to Speed up the Weight Optimization}

We first establish the proposed surrogate model and then introduce an efficient method for optimizing the weight using this model.

\subsubsection{Multilayer Perceptron-Based Surrogate Model}

\begin{figure}[!t]
   \centering
   \includegraphics[width=0.47\textwidth]{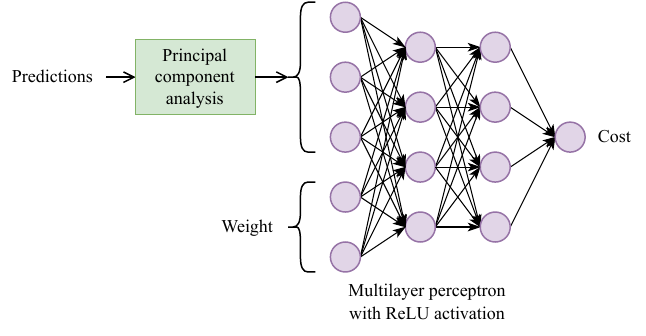}
   \caption{The proposed surrogate model.}
   \label{fig:surrogate}
\end{figure}

Fig.~\ref{fig:surrogate} illustrates the proposed surrogate model, which employs a multilayer perceptron (MLP) neural network to capture the mapping from predictions and weight to performance (i.e., the evaluated cost). The predictions $( \hat{\mathbf{U}}^{(m)}; m \in \mathcal{M})$ are high-dimensional, so we initially apply principal component analysis (PCA) \cite{abdi2010principal} for dimensionality reduction. The resulting components are then incorporated as inputs to the neural network. Since the weight $w = (w^{(m)}; m \in \mathcal{M})$ satisfies $\sum_{m \in \mathcal{M}} w^{(m)} = 1$, we input the first $|\mathcal{M}| - 1$ weight components into the neural network, where $|\mathcal{M}|$ denotes the number of elements in $\mathcal{M}$. 
The inputs from the predictions and the weight are concatenated and scaled before being processed by the neural network. 
The rectified linear unit (ReLU) activation function is employed, which outputs $v := \max \{0, s\}$ for an input scalar $s$. The neural network's output represents the evaluated cost. 
The neural network is trained offline, and once trained, the surrogate model can forecast the evaluated cost based on the provided predictions and weight data.

\subsubsection{MILP-Based Weight Optimization}

We model the proposed surrogate model using mixed-integer linear constraints. These linear constraints represent PCA, scaling, and linear combinations. The ReLU activation function $v = \max \{0, s\}$ can be equivalently expressed through the following set of constraints \cite{anderson2020strong}:
\begin{subequations}
\label{eq:ReLU}
\begin{align}
    \label{eq:ReLU-a}
    & 0 \leq v \leq M z, \\
    \label{eq:ReLU-b}
    & s \leq v \leq s + M (1 - z),
\end{align}
\end{subequations}
where $z$ is an auxiliary binary variable and $M > 0$ is a sufficiently large constant.
When $z = 0$, \eqref{eq:ReLU-a} forces $v = 0$ and \eqref{eq:ReLU-b} implies $v \geq s$; when $z = 1$, \eqref{eq:ReLU-b} shows $v = s$ and by \eqref{eq:ReLU-a} we have $v \geq 0$. 
Therefore, the ReLU activation function is equivalently modeled by the linear constraints in (4) with the binary variable $z$. The surrogate model is therefore expressed through these mixed-integer linear constraints:
\begin{subequations}
\label{eq:surrogate}
\begin{align}
\label{eq:surrogate-1}
& \hat{\mathbf{U}} = \sum_{m \in \mathcal{M}} w^{(m)} \hat{\mathbf{U}}^{(m)}, \\
\label{eq:surrogate-2}
& v_{0 j} = W_{0 j} P_j \hat{\mathbf{U}} + B_{0 j}, j = 1, 2, \dots, N_P, \\
\label{eq:surrogate-3}
& v_{0 (m + N_P)} = W_{0 (m + N_P)} w^{(m)} + B_{0 (m + N_P)}, \nonumber \\
& m = 1, 2, \dots, |\mathcal{M}| - 1, \\
\label{eq:surrogate-4}
& s_{i j} = W_{i j} v_{i - 1} + B_{i j}, i = 1, 2, \dots, N_L, j = 1, 2, \dots, N_{L i}, \\
\label{eq:surrogate-5}
& 0 \leq v_{i j} \leq M z_{i j}, i = 1, 2, \dots, N_L, j = 1, 2, \dots, N_{L i}, \\
\label{eq:surrogate-6}
& s_{i j} \leq v_{i j} \leq s_{i j} + M(1 - z_{i j}), \nonumber \\
& i = 1, 2, \dots, N_L, j = 1, 2, \dots, N_{L i}, \\
\label{eq:surrogate-7}
& z_{i j} \in \{0, 1\}, i = 1, 2, \dots, N_L, j = 1, 2, \dots, N_{L i}, \\
\label{eq:surrogate-8}
& I = W_{N_L + 1} v_{N_L} + B_{N_L + 1}. 
\end{align}
\end{subequations}
The combined prediction is represented in \eqref{eq:surrogate-1}. The PCA process is modeled in \eqref{eq:surrogate-2} with PCA parameters $P_1, P_2, \dots, P_{N_P}$, scaling parameters $W_{0 1}, W_{0 2}, \dots, W_{0 N_P}$, and bias parameters $B_{0 1}, B_{0 2}, \dots, B_{0 N_P}$. Similarly, the first $|\mathcal{M}| - 1$ components of the weight $w$ are scaled in \eqref{eq:surrogate-3}. The values from \eqref{eq:surrogate-2} and \eqref{eq:surrogate-3} are concatenated to form the inputs of the neural network, denoted as $v_0 = (v_{0 1}, v_{0 2}, \dots, v_{0 (N_P + |\mathcal{M}| - 1)})$. Constraint \eqref{eq:surrogate-4} models the relationship between the output $v_{i - 1}$ of the $(i-1)$-th layer and the input $s_i$ of $i$-th layer's ReLU function, where $W_{i j}$ and $B_{i j}$ are weight and bias parameters. The neural network consists of $N_L$ layers, with the $i$-th layer containing $N_{L i}$ units. Constraints \eqref{eq:surrogate-5} and \eqref{eq:surrogate-6} model the ReLU activation function, utilizing the binary variable $z_{i j}$ defined in \eqref{eq:surrogate-7}. Constraint \eqref{eq:surrogate-8} specifies the output $I$ of the surrogate model.

To efficiently optimize the weight, we solve the MILP problem \eqref{eq:weight-optimization} to minimize the output of the neural network:
\begin{align}
\min_{w, \hat{\mathbf{U}}, s, v, z, I}~ & I \nonumber \\
\label{eq:weight-optimization}
\mbox{s.t.}~ & (5), w \geq 0, \sum_{m \in \mathcal{M}} w^{(m)} = 1. 
\end{align}
The optimal weight value is then utilized in data-driven RO to derive the final strategy.

We use the proposed surrogate model based on MLP with ReLU activation to describe how the predictions and weight influence the performance of the optimized strategy. MLP with a non-polynomial activation function can approximate any continuous function to any desired degree of accuracy \cite{leshno1993multilayer}. Additionally, ref. \cite{lu2017expressive} examined the width and depth needed for MLP with ReLU activation to serve as universal approximators. These theoretical results demonstrate the underlying rationale for the proposed MLP-based surrogate model.

\subsection{Overview of the Proposed Framework}

The steps of the proposed integrated forecasting and optimization framework are summarized as follows:

Step 1 (offline training): Train the surrogate model using historical data on uncertainty. Details are provided in Algorithm~\ref{alg:surrogate}.

\begin{algorithm}[!t]
\label{alg:surrogate}
\DontPrintSemicolon 
  \SetAlgoLined
  \KwIn{Uncertainty samples $u_1, u_2, \dots, u_{N_U}$; weight samples $w_1, w_2, \dots, w_{N_W}$; dimension $N_P$ of the PCA output; number of layers $N_L$ in the MLP; number of units $N_{L i}$ in layer $i$, where $i = 1, 2, \dots, N_L$.}
  \KwOut{PCA parameters $P_1, P_2, \dots, P_{N_P}$; MLP weight parameter $W_{i j}$ and bias parameter $B_{i j}$ for $i = 0, 1, \dots, N_L$ and $j = 1, 2, \dots, N_{L i}$.}
  Compute PCA parameters $P_1, P_2, \dots, P_{N_P}$ using uncertainty samples $u_1, u_2, \dots, u_{N_U}$.\;
  \For{$n_U = 1$ \textbf{to} $N_U$}{
    Compute PCA outputs $d_{n_U 1} = P_1 u_{n_U}$, $d_{n_U 2} = P_2 u_{n_U}$, \dots, $d_{n_U N_P} = P_{N_P} u_{n_U}$.\;
    \For{$n_W = 1$ \textbf{to} $N_W$}{
    Compute $\hat{\mathbf{U}}(w_{n_W}) = \sum_{m \in \mathcal{M}} w_{n_W}^{(m)} \hat{\mathbf{U}}^{(m)}$ and $\mathbf{E}^{(1:N)}(w_{n_W}) = \mathbf{U}^{(1:N)} - \hat{\mathbf{U}}^{(1:N)}(w_{n_W})$.\;
    Based on $\hat{\mathbf{U}}(w_{n_W})$ and $\mathbf{E}^{(1:N)}(w_{n_W})$, apply Algorithm 2 to obtain the uncertainty set $\mathcal{U}(w_{n_W})$ and the optimal solution $x_\mathcal{U}^*(w_{n_W})$.\;
    Using $x_\mathcal{U}^*(w_{n_W})$, compute the evaluated cost $I(u_{n_U}, w_{n_W})$ as described in Section~\ref{sec:framework-integration}.\;
    }
  }
  Use the input data $(d_{n_U 1}, d_{n_U 2}, \dots, d_{n_U N_P}, w_{n_W})$ and the output data $I(u_{n_U}, w_{n_W})$ for $n_U = 1, 2, \dots, N_U$ and $n_W = 1, 2, \dots, N_W$ to train the MLP with ReLU activation, and obtain the parameters $W_{i j}$ and $B_{i j}$ for $i = 1, 2, \dots, N_L$ and $j = 1, 2, \dots, N_{L i}$.\;
  \caption{Surrogate model training}
\end{algorithm}

Step 2 (daily optimization for UC): Solve the MILP problem \eqref{eq:weight-optimization} to obtain the optimal weight $w$ and the combined prediction $\hat{\mathbf{U}}$. Solve the data-driven RO problem \eqref{eq:RUC-original} using Algorithm~\ref{alg:RUC} and obtain the final strategy.

Compared to existing predict-and-optimize frameworks in power system applications, the proposed method has several advantages: First, instead of adjusting the parameters of forecasting models, the proposed framework optimizes the weight $w \in \mathcal{W}$ assigned to different forecasting methods. Thus, the prediction used for optimization must be a convex combination of the original predictions, which prevents severe deviations and enhances the reliability of the induced strategy. 
Second, we introduce a MLP-based surrogate model to represent the relationship between the weight and evaluated cost. This allows for rapid weight optimization via MILP, which significantly reduces computational burdens compared to existing frameworks that rely on forward and backward propagation or other iterative algorithms. Third, the proposed method incorporates data-driven RO in the optimization process, which provides statistical guarantees within the integrated forecasting and optimization framework. This approach ensures the out-of-sample performance of the optimized strategy from a theoretical perspective.

\section{Data-Driven Two-Stage Robust Unit Commitment with Statistical Guarantees}
\label{sec:RUC}

In this section, we first establish the robust UC problem and then construct the uncertainty set, which is later reconstructed using the problem information. The solution algorithm is introduced in the end.

\subsection{Robust Unit Commitment Formulation}
\label{sec:RUC-problem}

As introduced in Section~\ref{sec:framework-optimization}, the robust UC problem has the compact form in \eqref{eq:RUC-original}. Now we specify the components of this problem. The pre-dispatch variable $x = (\theta_{g t}, \theta_{g t}^\pm, p_{g t}, r_{g t}^\pm; g \in \mathcal{G}, t \in \mathcal{T})$, where $\mathcal{G}$ and $\mathcal{T}$ are the index sets of generators and periods, respectively. 
For generator $g$ in period $t$, the on, startup, and shutdown states are denoted by binary variables $\theta_{g t}$, $\theta_{g t}^+$, and $\theta_{g t}^-$, respectively; the day-ahead scheduled power output is denoted by $p_{g t}$; the upward and downward reserve power values are denoted by $r_{g t}^+$ and $r_{g t}^-$, respectively.
The pre-dispatch cost is given by
\begin{align}
    \label{eq:pre-cost}
    f(x) = \sum_{t \in \mathcal{T}} \sum_{g \in \mathcal{G}} \left( o_g^+ \theta_{g t}^+ + o_g^- \theta_{g t}^- + \rho_g p_{g t} + \gamma_g^+ r_{g t}^+ + \gamma_g^- r_{g t}^- \right),
\end{align}
where $o_g^\pm$, $\rho_g$, and $\gamma_g^\pm$ are cost coefficients. The pre-dispatch feasible region is defined as follows:
\begin{subequations}
\label{eq:X}
    \begin{align}
        & \mathcal{X} = \left\{ x = (\theta_{g t}, \theta_{g t}^\pm, p_{g t}, r_{g t}^\pm; g \in \mathcal{G}, t \in \mathcal{T}) ~\right| \nonumber \\
        \label{eq:X-a}
        & \sum_{g \in \mathcal{G}} p_{g t} = \sum_{i \in \mathcal{I}} \hat{u}_{i t}, \forall t \in \mathcal{T}, \\
        \label{eq:X-b}
        & - S_l \leq \sum_{g \in \mathcal{G}} \pi_{g l} p_{g t} - \sum_{i \in \mathcal{I}} \pi_{i l} \hat{u}_{i t} \leq S_l, \forall l \in \mathcal{L}, \forall t \in \mathcal{T}, \\
        \label{eq:X-c}
        & \theta_{g t},\theta_{g t}^+,\theta_{g t}^- \in \{0,1\}, \forall g \in \mathcal{G}, \forall t \in \mathcal{T}, \\
        \label{eq:X-d}
        & \sum_{\tau = t}^{t + T_g^+ - 1} \theta_{g \tau} \geq T_g^+ \theta_{gt}^+, 1 \leq t \leq T - T_g^+ + 1, \forall g \in \mathcal{G}, \\
        \label{eq:X-e}
        & \sum_{\tau = t}^T (\theta_{g \tau} - \theta_{g t}^+) \geq 0, T - T_g^+ + 2 \leq t \leq T, \forall g \in \mathcal{G}, \\
        \label{eq:X-f}
        & \sum_{\tau = t}^{t + T_g^- - 1} (1 - \theta_{g \tau}) \geq T_g^- \theta_{g t}^-, 1 \leq t \leq T - T_g^- + 1, \forall g \in \mathcal{G}, \\
        \label{eq:X-g}
        & \sum_{\tau = t}^T (1 - \theta_{g \tau} - \theta_{g t}^-) \geq 0, T - T_g^- + 2 \leq t \leq T, \forall g \in \mathcal{G}, \\
        \label{eq:X-h}
        & \theta_{g t} - \theta_{g (t-1)} = \theta_{g t}^+ - \theta_{g t}^-, \forall g \in \mathcal{G}, \forall t \in \mathcal{T}, \\
        \label{eq:X-i}
        & \theta_{g t}^+ + \theta_{g t}^- \leq 1, \forall g \in \mathcal{G}, \forall t \in \mathcal{T}, \\
        \label{eq:X-j}
        & 0 \leq r_{g t}^+ \le R_g^+ \theta_{g t}, 0 \leq r_{g t}^- \leq R_g^- \theta_{g t}, \forall g \in \mathcal{G}, \forall t \in \mathcal{T}, \\
        \label{eq:X-k}
        & \underline{P}_g \theta_{g t} + r_{g t}^- \leq p_{g t} \leq \overline{P}_g \theta_{g t} - r_{g t}^+, \forall g \in \mathcal{G}, \forall t \in \mathcal{T}, \\
        & (p_{g t} + r_{g t}^+) - (p_{g (t-1)} - r_{g (t-1)}^-) \leq K_g^+ \theta_{g (t-1)} \nonumber\\ 
        \label{eq:X-l}
        &  + K_g^U \theta_{g t}^+, \forall g \in \mathcal{G}, 2 \leq t \leq T, \\
        & - (p_{g t}-r_{g t}^-) + (p_{g (t-1)} + r_{g (t-1)}^+) \leq K_g^- \theta_{g t} \nonumber\\
        \label{eq:X-m}
        & \left. + K_g^D \theta_{g t}^-, \forall g \in \mathcal{G}, 2 \leq t \leq T \right\}.
    \end{align}
\end{subequations}
The power balance of the transmission network is stipulated in \eqref{eq:X-a}, where $\hat{u}_{i t}$ is the day-ahead prediction of the load power. In \eqref{eq:X-b}, $\mathcal{L}$ is the index set of transmission lines; $S_l$ is the power capacity of line $l$; 
parameters $\pi_{g l}$ and $\pi_{i l}$ are the power transfer distribution factors. 
Thus, \eqref{eq:X-b} bounds the line flow in the DC power flow model. Binary variables are set in \eqref{eq:X-c}. Constraints \eqref{eq:X-d}-\eqref{eq:X-g} are for the minimum up time $T_g^+$ and minimum down time $T_g^-$ of generators \cite{arroyo2000optimal}. The generator state change is modeled in \eqref{eq:X-h} and the simultaneous startup and shutdown of a generator is prohibited in \eqref{eq:X-i}. Constraint \eqref{eq:X-j} contains bounds $R_g^\pm$ for the upward and downward reserve power. The minimum and maximum outputs of generators, i.e., $\underline{P}_g$ and $\overline{P}_g$ for generator $g$, are stipulated in \eqref{eq:X-k}. For generator $g$, $K_g^+$ and $K_g^-$ are the maximum upward and downward ramp values in a period; $K_g^U$ and $K_g^D$ are the maximum increase and decrease of power if the generator startups or shutdowns in a period. Thus, constraints \eqref{eq:X-l} and \eqref{eq:X-m} are the bounds for the change in generator power over a period.

The construction of the uncertainty set $\mathcal{U}$ is deferred and will be addressed in Sections~\ref{sec:RUC-uncertainty-1} and \ref{sec:RUC-uncertainty-2}. In the re-dispatch stage, the variable is $y = (p_{g t}^\pm; g \in \mathcal{G}, t \in \mathcal{T})$, where $p_{g t}^+$ and $p_{g t}^-$ are the upward and downward power adjustments of generator $g$ in period $t$, respectively. The re-dispatch cost function is the total power adjustment cost, i.e.,
\begin{align}
    \label{eq:re-cost}
    h(y) = \sum\nolimits_{t \in \mathcal{T}} \sum\nolimits_{g \in \mathcal{G}} \left( \rho_g^+ p_{g t}^+ + \rho_g^- p_{g t}^- \right), 
\end{align}
where $\rho_g^+$ and $\rho_g^-$ are cost coefficients. The re-dispatch feasible region $\mathcal{Y}(x, u)$ depends on the pre-dispatch decision $x$ and the realization $u = (u_{i t}; i \in \mathcal{I}, t \in \mathcal{T})$ of uncertain load:
\begin{subequations}
    \label{eq:Y}
    \begin{align}
        & \mathcal{Y}(x, u)= \left\{ y = (p_{g t}^\pm; g \in \mathcal{G}, t \in \mathcal{T}) ~\right| \nonumber\\
        \label{eq:Y-a}
        & \sum\nolimits_{g \in \mathcal{G}} (p_{g t}+ p_{g t}^+ - p_{g t}^-) = \sum\nolimits_{i \in \mathcal{I}} u_{i t}, \forall t \in \mathcal{T}, \\
        & - S_l \leq \sum\nolimits_{g \in \mathcal{G}} \pi_{g l} (p_{g t}+ p_{g t}^+ - p_{g t}^-) - \sum\nolimits_{i \in \mathcal{I}} \pi_{i l} u_{i t} \nonumber \\
        \label{eq:Y-b}
        & \leq S_l, \forall l \in \mathcal{L}, \forall t \in \mathcal{T}, \\
        \label{eq:Y-c}
        & \left. 0 \leq p_{g t}^+ \leq r_{g t}^+, 0 \leq p_{g t}^- \leq r_{g t}^-, \forall g \in \mathcal{G}, \forall t \in \mathcal{T} \right\}.
    \end{align}
\end{subequations}
In \eqref{eq:Y}, constraint \eqref{eq:Y-a} is for the power balance of the network. The transmission line flow constraints are in \eqref{eq:Y-b}. The power adjustments of the generators are bounded by the reserve power in \eqref{eq:Y-c}.

It is crucial to consider ramp rate constraints to ensure the feasibility of the UC strategy \cite{moreira2024role}. Lemma~\ref{lemma:ramp} shows that the re-dispatch feasible region $\mathcal{Y}(x,u)$ adheres the ramp rate constraints, provided that $x \in \mathcal{X}$. The proof of Lemma~\ref{lemma:ramp} is in Appendix~\ref{appendix:lemma-ramp}.

\begin{lemma}
\label{lemma:ramp}
    Assume $x = (\theta_{g t}, \theta_{g t}^\pm, p_{g t}, r_{g t}^\pm; g \in \mathcal{G}, t \in \mathcal{T}) \in \mathcal{X}$ and $y = (p_{g t}^\pm; g \in \mathcal{G}, t \in \mathcal{T}) \in \mathcal{Y}(x,u)$. Then the following ramp rate constraints hold for the re-dispatch solution:
    \begin{subequations}
    \begin{align}
        & (p_{g t}+ p_{g t}^+ - p_{g t}^-) - (p_{g (t-1)}+ p_{g (t-1)}^+ - p_{g (t-1)}^-) \nonumber \\ 
        \label{eq:ramp-re-upper}
        & \leq K_g^+ \theta_{g (t-1)} + K_g^U \theta_{g t}^+, \forall g \in \mathcal{G}, 2 \leq t \leq T, \\
        & - (p_{g t}+ p_{g t}^+ - p_{g t}^-) + (p_{g (t-1)}+ p_{g (t-1)}^+ - p_{g (t-1)}^-) \nonumber\\
        \label{eq:ramp-re-lower}
        & \leq K_g^- \theta_{g t} + K_g^D \theta_{g t}^-, \forall g \in \mathcal{G}, 2 \leq t \leq T.
    \end{align}
    \end{subequations}
\end{lemma}

According to the formulations of $f(x)$ and $h(y)$ in \eqref{eq:pre-cost} and \eqref{eq:re-cost}, the two functions are linear. Equation \eqref{eq:Y} shows that the constraints that define $\mathcal{Y}(x, u)$ are linear in $x$ and $u$. Therefore, the compact form of the robust UC problem can be further written as follows:
\begin{align}
    \label{eq:RUC-linear}
    \min\nolimits_{x \in \mathcal{X}} \left\{ C^\top x + \max\nolimits_{u \in \mathcal{U}} \min\nolimits_{y: A y \geq B x + D u + E} F^\top y \right\},
\end{align}
where $A$, $B$, $C$, $D$, $E$, and $F$ are coefficient matrices and vectors.

\subsection{Data-Driven Uncertainty Set and Statistical Guarantees}
\label{sec:RUC-uncertainty-1}

The possible realization values of the uncertain load $u$ constitute a set $\mathcal{U}_0 = \{ u ~|~ \underline{U} \leq u \leq \overline{U} \}$, where $\underline{U}$ and $\overline{U}$ are the lower and upper bounds for $u$. 
$\mathcal{U}_0$ contains all the possible values of $u$, so $\Pr [u \in \mathcal{U}_0] = 1$.
$\mathcal{U}_0$ can be constructed using the capacity parameters of the facilities. For example, if bus $i$ connects a wind farm with a power capacity of $P_i^W$ and no loads, then the net load power at bus $i$ satisfies $-P_i^W \leq u_{i t} \leq 0$ for any $t \in \mathcal{T}$. The box defined by these bounds for all $i \in \mathcal{I}$ and $t \in \mathcal{T}$ establishes a choice of $\mathcal{U}_0$ such that $\Pr[u \in \mathcal{U}_0] = 1$. If the accurate $\mathcal{U}_0$ is hard to obtain, we may just let $\mathcal{U}_0$ to be the entire Euclidean space, and the proposed method will still function effectively.

According to the analysis in Section~\ref{sec:framework-optimization}, we need to construct an uncertainty set $\mathcal{U}_1 \subset \mathcal{U}_0$ such that $\Pr [u \in \mathcal{U}_1] \geq 1 - \varepsilon$, where the information we have include prediction $\hat{u}$ and historical forecast error $e_{1: N} := (e_1, e_2, \dots, e_N)$ in the training dataset ($e_n$ is the column vector reshaping $\mathbf{E}^{(n)}(w)$). We assume that the daily load forecast errors are i.i.d. continuous random variables. Our plan is first to construct a set $\mathcal{E}$ for the uncertain forecast error $e := u - \hat{u}$ so that $\Pr[e \in \mathcal{E}] \geq 1 - \varepsilon$, followed by letting
\begin{align}
    \label{eq:U1}
    \mathcal{U}_1 := \{ u \in \mathcal{U}_0 ~|~ u - \hat{u} \in \mathcal{E} \}. 
\end{align}
Clearly, such $\mathcal{U}_1$ satisfies $\mathcal{U}_1 \subset \mathcal{U}_0$ and $\Pr [u \in \mathcal{U}_1] \geq 1 - \varepsilon$, so the uncertainty set construction problem comes down to finding a set $\mathcal{E}$ such that $\Pr[e \in \mathcal{E}] \geq 1 - \varepsilon$, based on the historical data $e_{1: N}$. 

Since $e_{1: N}$ are random variables, the set $\mathcal{E}$ constructed using $e_{1: N}$ is also random, and so is the event $\Pr[e \in \mathcal{E}] \geq 1 - \varepsilon$. Therefore, instead of directly attempting $\Pr[e \in \mathcal{E}] \geq 1 - \varepsilon$, we consider the following statistical guarantee: 
\begin{align}
    \label{eq:guarantee-error}
    \mathbb{P}^N[\Pr[e \in \mathcal{E}] \geq 1 - \varepsilon] \geq 1 - \delta,
\end{align}
where $\mathbb{P}$ is the underlying distribution of $e_n, n = 1, 2, \dots, N$; $\mathbb{P}^N$ denotes the $N$ times distribution product of $\mathbb{P}$, which models the uncertainty of historical data and $\mathcal{E}$; $\delta$ is a probability tolerance parameter. Thus, equation \eqref{eq:guarantee-error} means that the probability of the random event $\Pr[e \in \mathcal{E}] \geq 1 - \varepsilon$ is at least $1 - \delta$. The parameters $\varepsilon$ and $\delta$ jointly control the conservative degree of the robust UC problem.

To achieve the statistical guarantee \eqref{eq:guarantee-error}, we adopt the data-driven uncertainty set proposed in \cite{hong2021learning} to establish $\mathcal{E}$, whose procedure and statistical guarantee are summarized in Theorem~\ref{thm:construction}. The idea is to pull two disjoint groups from the dataset. The first group determines the shape of $\mathcal{E}$, where an ellipsoid is established based on the sample mean and covariance to consider the correlation. The second group is for the size of $\mathcal{E}$, where the differences of the points from the ellipsoid center are measured and sorted to provide a threshold, and the independence of the two data groups forms the foundation of the statistical guarantees.

\begin{theorem}[Uncertainty set construction and statistical guarantee]
\label{thm:construction}
    Select two disjoint groups from $e_{1: N}$. Denote their realizations by $e_{1: N_1}^{(1)}$ and $e_{1: N_2}^{(2)}$, where $N_1 + N_2 \leq N$ and $N_2 \geq \log_{1 - \varepsilon} \delta$. Let
    \begin{subequations}
        \label{eq:sample-moment}
        \begin{align}
            & \mu := \frac{1}{N_1} \sum_{n = 1}^{N_1} e_n^{(1)}, \\
            & \Sigma := \frac{1}{N_1 - 1} \sum_{n = 1}^{N_1} \left( e_n^{(1)} - \mu \right) \left( e_n^{(1)} - \mu \right)^\top. 
        \end{align}
    \end{subequations}
    In other words, $\mu$ and $\Sigma$ are the sample mean and the sample covariance matrix in the first dataset $e_{1: N_1}^{(1)}$. Assume $\Sigma$ is invertible. Define function $a(e) := (e - \mu)^\top \Sigma^{-1} (e - \mu)$. Let
    \begin{align}
        \label{eq:n*}
        n_N^* := \min \left\{ n \in \mathbb{N} ~\middle|~ \sum\nolimits_{m = 0}^{n - 1} C_N^m (1 - \varepsilon)^m \varepsilon^{N - m} \geq 1 - \delta \right\}, 
    \end{align}
    where $C_N^m$ denotes the binomial coefficient of $N$ choose $m$. Let $\alpha$ be the $n_{N_2}^*$-th smallest value in $a(e_1^{(2)}), a(e_2^{(2)}), \dots, a(e_{N_2}^{(2)})$. Define the set $\mathcal{E}$ as
    \begin{align}
        \label{eq:E}
        \mathcal{E} := \left\{ e ~\middle|~ (e - \mu)^\top \Sigma^{-1} (e - \mu) \leq \alpha \right\}.
    \end{align}
    Then the statistical guarantee \eqref{eq:guarantee-error} holds. Moreover, if $\mathcal{U} = \mathcal{U}_1$ defined in \eqref{eq:U1} is used in the robust UC problem \eqref{eq:RUC-original} and $x_0 := x_{\mathcal{U}_1}^*$ is an optimal solution, then 
    \begin{align}
        \label{eq:guarantee-construction}
        \mathbb{P}^N[O \leq O_{x_0} \leq O_{\mathcal{U}_1}] \geq 1 - \delta.
    \end{align}
\end{theorem}

Theorem~\ref{thm:construction} gives a method to generate the uncertainty set $\mathcal{U}_1$ and approximate problem \eqref{eq:UC-original} with confidence $1 - \delta$. The proof of Theorem~\ref{thm:construction} can be found in Appendix~\ref{appendix:thm-construction}.

Net loads often exhibit strong spatial and temporal correlations, making uncertainty correlation modeling essential for decision-making under multi-dimensional uncertainty \cite{velloso2019two}. The ellipsoidal uncertainty set $\mathcal{U}_1$ of $u = (u_{i t}; i \in \mathcal{I}, t \in \mathcal{T})$ is able to consider the correlation between $u_{i t}$ and $u_{j s}$ for any $i, j \in \mathcal{I}$ and $t, s \in \mathcal{T}$ through the sample covariance matrix $\Sigma$ of the forecast errors \cite{golestaneh2018ellipsoidal}.

\subsection{Uncertainty Set Reconstruction}
\label{sec:RUC-uncertainty-2}

The uncertainty set $\mathcal{U}_1$ obtained in Theorem~\ref{thm:construction} guarantees robustness with confidence $1 - \delta$, but it can be conservative. Our goal in this subsection is to mitigate the conservativeness by reconstructing the uncertainty set. The construction of $\mathcal{U}_1$ uses only the prediction and historical forecast error but does not involve any information on the UC problem. In the following, we integrate the problem information into the reconstruction. 

\begin{figure}[!t]
   \centering
   \includegraphics[width=0.28\textwidth, trim= 18mm 30mm 37mm 30mm, clip]{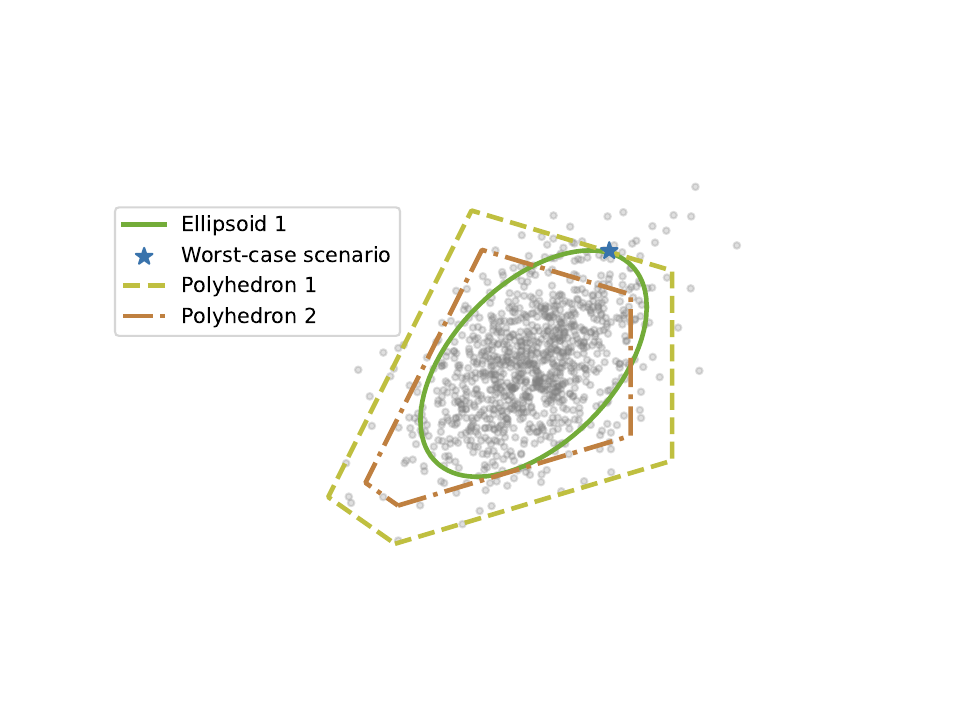}
   \caption{The basic idea of uncertainty set reconstruction.}
   \label{fig:sketch}
\end{figure}

We illustrate the basic idea in Fig.~\ref{fig:sketch}. 
Ellipsoid 1 represents $\mathcal{U}_1$. The worst-case scenario in ellipsoid 1 is marked as a blue star, which reaches the highest cost, with nearby scenarios also incurring high costs. 
Considering the range where costs do not exceed that of the blue star, we obtain polyhedron 1 in Fig.~\ref{fig:sketch}. This range forms a polyhedron due to the linear structure in the UC problem. Polyhedron 1 includes ellipsoid 1, and the blue star remains the worst-case scenario within polyhedron 1. Furthermore, polyhedron 1 may contain more data points than ellipsoid 1, which can lead to conservativeness. This suggests that polyhedron 1 can be shrunk into polyhedron 2 according to the desired probability guarantee thresholds. 
In contrast to ellipsoid 1, polyhedron 2 excludes the high-cost scenarios near the blue star, while including additional low-cost scenarios to maintain the statistical guarantees. Consequently, the worst-case scenario in polyhedron 2 will not exceed the blue star, resulting in a less conservative strategy.

Now we explain the reconstruction formally. Recall from Lemma~\ref{lemma:RUC} that when $(x^*, \eta^*)$ is optimal for the chance-constrained problem \eqref{eq:UC-original}, the uncertainty set $\mathcal{U}^*$ defined in \eqref{eq:U*} can equivalently transform problem \eqref{eq:UC-original} into \eqref{eq:RUC-original}. Ideally, $\mathcal{U}^*$ is used as the uncertainty set to eliminate conservativeness. However, it is impractical because $(x^*, \eta^*)$ is unknown. To this end, $\mathcal{U}^*$ is approximated using the data we have. Suppose $x_0$ is a feasible solution, then we can estimate its performance $O_{x_0}$ in a historical dataset and construct an uncertainty set using $(x, \eta) = (x_0, O_{x_0} - f(x_0))$. This idea is refined to maintain the statistical guarantees in Theorem~\ref{thm:reconstruction}.

\begin{theorem}[Uncertainty set reconstruction]
    \label{thm:reconstruction}
    Suppose $e_{1: N_3}^{(3)}$ is a subgroup of $e_{1: N}$ and $N_3 \geq \log_{1 - \varepsilon} \delta$. Assume that the solution $x_0$ is independent of $e_{1: N_3}^{(3)}$. Define function $b(u) := \min_{y \in \mathcal{Y}(x_0, u)} h(y)$. Define $n_{N_3}^*$ according to \eqref{eq:n*}. Let $\beta$ be the $n_{N_3}^*$-th smallest value in $b(\hat{u} + e_1^{(3)}), b(\hat{u} + e_2^{(3)}), \dots, b(\hat{u} + e_{N_3}^{(3)})$. Define the uncertainty set
    \begin{align}
        \label{eq:U2}
        \mathcal{U}_2 := \left\{ u \in \mathcal{U}_0 ~\middle|~ \min\nolimits_{y \in \mathcal{Y}(x_0, u)} h(y) \leq \beta \right\}.
    \end{align}
    Let $x_1 := x_{\mathcal{U}_2}^*$ be optimal in problem \eqref{eq:RUC-original} with $\mathcal{U} = \mathcal{U}_2$. Then the following statistical guarantees hold:
    \begin{subequations}
        \label{eq:guarantee-reconstruction}
        \begin{align}
            \label{eq:guarantee-reconstruction-1}
            & \mathbb{P}^N[\Pr[u \in \mathcal{U}_2] \geq 1 - \varepsilon] \geq 1 - \delta, \\
            \label{eq:guarantee-reconstruction-2}
            & \mathbb{P}^N[O \leq O_{x_1} \leq O_{\mathcal{U}_2} \leq f(x_0) + \beta] \geq 1 - \delta.
        \end{align}
    \end{subequations}
\end{theorem}

In Theorem~\ref{thm:reconstruction}, the shape of the new uncertainty set is formed based on both the old solution and the UC problem, whereas the size $\beta$ is determined by the evaluations on an independent dataset so that the statistical guarantees remain valid. The conclusion indicates that the reconstruction of the uncertainty set will lead to a new solution probably better than the old one, represented by $O_{x_1} \leq f(x_0) + \beta$. The proof of Theorem~\ref{thm:reconstruction} can be found in Appendix~\ref{appendix:thm-reconstruction}.

The polyhedral uncertainty set $\mathcal{U}_2$ also takes into account the spatial and temporal correlations of net loads. However, unlike the ellipsoidal uncertainty set $\mathcal{U}_1$, which uses the sample covariance matrix of historical data, $\mathcal{U}_2$ describes correlations in a problem-driven way, through linear constraints derived from the UC problem. As will be elaborated in Fig.~\ref{fig:projection} in Section~\ref{sec:case}, $\mathcal{U}_2$ is not box-shaped, demonstrating that it incorporates correlations across different dimensions.

\subsection{Solution Algorithm}
\label{sec:RUC-solution}

$\mathcal{U}_0$ is polyhedral. 
According to \eqref{eq:U1} and \eqref{eq:E},
\begin{align}
    \mathcal{U}_1 = \{ u \in \mathcal{U}_0 ~|~ (u - \hat{u} - \mu)^\top \Sigma^{-1} (u - \hat{u} - \mu) \leq \alpha \}, \nonumber
\end{align}
is the intersection of an ellipsoid and a polyhedron. Using the compact form of the robust UC problem in \eqref{eq:RUC-linear}, $\mathcal{U}_2$ in \eqref{eq:U2} can be further written as follows:
\begin{align}
    \mathcal{U}_2 & = \{ u \in \mathcal{U}_0 ~|~ \exists y \in \mathcal{Y}(x_0, u), ~\mbox{s.t.}~ h(y) \leq \beta \} \nonumber \\
    & = \{ u \in \mathcal{U}_0 ~|~ \exists y, ~\mbox{s.t.}~ A y \geq B x_0 + D u + E, F^\top y \leq \beta \}, \nonumber
\end{align}
which shows that $\mathcal{U}_2$ is polyhedral. Problem \eqref{eq:RUC-linear} with these uncertainty sets can be effectively solved by the column-and-constraint generation (C\&CG) algorithm \cite{zeng2013solving}, which is omitted here for the sake of conciseness.

The solution procedure for the robust UC problem considering uncertainty set reconstruction is summarized in Algorithm~\ref{alg:RUC}. The historical dataset is divided into two groups. The first group forms the first ellipsoidal uncertainty set and leads to a solution $x_0$. Then the reconstruction procedure in Theorem~\ref{thm:reconstruction} is performed to obtain an improved uncertainty set, resulting in the final solution $x_1$. Since the two datasets are independent, the statistical guarantees in \eqref{eq:guarantee-reconstruction} are maintained.

\begin{algorithm}[!t]
\label{alg:RUC}
\DontPrintSemicolon 
  \SetAlgoLined
  \KwIn{Parameters of \eqref{eq:RUC-linear}; $\varepsilon$; $\delta$; $\mathcal{U}_0$; $\hat{u}$; $e_{1: N}$; $N_2$.}
  \KwOut{UC strategy $x_1$.}
  Let $N_1 \leftarrow N - N_2$.\;
  Divide $e_{1: N}$ into $e_{1: N_1}^{(1)}$ and $e_{1: N_2}^{(2)}$.\;
  Calculate $\mu$ and $\Sigma$ according to \eqref{eq:sample-moment}.\;
  Let $\alpha \leftarrow \max \{ (e_n^{(1)} - \mu)^\top \Sigma^{-1} (e_n^{(1)} - \mu) ~|~ n = 1, 2, \dots, N_1 \}$.\;
  Let $\mathcal{U}_1' \leftarrow \{ u \in \mathcal{U}_0 | (u - \hat{u} - \mu)^\top \Sigma^{-1} (u - \hat{u} - \mu) \leq \alpha \}$.\;
  Solve problem \eqref{eq:RUC-linear} with $\mathcal{U} = \mathcal{U}_1'$ using the C\&CG algorithm and obtain the optimal solution $x_0$.\;
  Let $b_n \leftarrow \min_{y: A y \geq B x_0 + D (\hat{u} + e_n^{(2)}) + E} F^\top y$, for $n = 1, 2, \dots, N_2$.\;
  Arrange $b_n, n = 1, 2, \dots, N_2$ from small to large and get $b_n', n = 1, 2, \dots, N_2$.\;
  Let $n^* \leftarrow \min \{ n | \sum_{m = 0}^{n - 1} C_{N_2}^m (1 - \varepsilon)^m \varepsilon^{N_2 - m} \geq 1 - \delta \}$.\;
  Let $\beta \leftarrow b_{n^*}'$.\;
  Let $\mathcal{U}_2 \leftarrow \{ u \in \mathcal{U}_0 | \exists y, ~\mbox{s.t.}~ A y \geq B x_0 + D u + E, F^\top y \leq \beta \}$.\;
  Solve problem \eqref{eq:RUC-linear} with $\mathcal{U} = \mathcal{U}_2$ using the C\&CG algorithm and return the optimal solution $x_1$.\;
  \caption{Solution of robust unit commitment}
\end{algorithm}

\emph{Remark}: The proposed method does not require explicit consideration of the robust feasibility requirement and can effectively address the infeasible re-dispatch cases. For a fixed first-stage variable $x' \in \mathcal{X}$, if the re-dispatch problem is infeasible for a given uncertainty realization $u' \in \mathcal{U}$, then $\mathcal{Y}(x', u') = \varnothing$ and the second-stage cost
\begin{align}
    \max_{u \in \mathcal{U}} \min_{y \in \mathcal{Y}(x', u)} h(y) \geq \min_{y \in \mathcal{Y}(x', u')} h(y) = + \infty. \nonumber
\end{align}
Since the two-stage RO problem \eqref{eq:RUC-original} minimizes the sum of the first- and second-stage costs, the optimal solution $x^*$ must be robust feasible, i.e., $\mathcal{Y}(x^*, u) \neq \varnothing$ for all $u \in \mathcal{U}$, provided that the optimal value is finite. Therefore, the robust feasibility requirement naturally holds for the optimal solution, effectively excluding scenarios with infeasible re-dispatch.

\section{Case Studies}
\label{sec:case}

This section examines the proposed robust UC method using modified IEEE 30-bus and 118-bus systems. All experiments are carried out on a laptop with an Intel i7-12700H processor and 16 GB RAM. The neural network is established and trained by PyTorch 2.1.2. The MILP problems in the C\&CG algorithm are solved by Gurobi 11.0.2. In the following, we first introduce the prediction data. The performance of the proposed method is then investigated in the modified IEEE 30-bus system, where different methods are compared and sensitivity analysis is performed. Finally, the modified IEEE 118-bus system is used to demonstrate the scalability of the proposed method.

\subsection{Prediction Data}

The hourly load data in one and a half years are extracted from the dataset in \cite{grabner2023global}, based on real data in Ireland \cite{commission2012cer}. The hourly wind power data are generated according to historical weather data \cite{staffell2016using}. The load and wind data are aligned according to the date information. 
The dataset is divided into samples for training, validation, and testing. 

We predict uncertainty using the following forecasting methods:
\begin{itemize}
    \item M1: Forecast nodal power using local data and bidirectional long short-term memory (BiLSTM) neural network \cite{siami2019performance}.
    \item M2: Apply federated learning between nodes and forecast using the adapted global BiLSTM model \cite{fekri2022distributed}.
    \item M3: Forecast nodal power using BiLSTM networks and nodal power subprofiles \cite{wang2018ensemble}.
    \item C1: Combine the predictions of M1, M2, and M3 by a weight $w$ to minimize the mean square error (MSE).
    \item C2: Combine the predictions of M1, M2, and M3 by a weight $w$ according to the method proposed in Section~\ref{sec:framework}.
\end{itemize}
The loss function used in the training is MSE. We adopt root mean square error (RMSE) and mean absolute error (MAE) to measure the test forecast errors.

The average forecast errors under different methods are shown in TABLE~\ref{tab:load}. By combining three forecasting methods and optimizing the weight to minimize the MSE, C1 has the lowest forecast errors, where its RMSE is 11.1\%, 7.3\%, and 6.8\% lower than M1, M2, and M3, respectively. This verifies that the combination of predictions is effective in improving forecast accuracy by leveraging the potential of different data sources and forecasting methods. The forecast errors of C2 are larger than those of C1 because instead of minimizing the MSE, the weight in C2 is chosen to minimize the UC cost.

\begin{table}[!t]
\renewcommand{\arraystretch}{1.3}
\caption{Average Forecast Errors of Different Methods}
\label{tab:load}
\centering
\begin{tabular}{ccc}
\hline
Method & RMSE & MAE \\
\hline
M1 & 84.39 & 54.64 \\
M2 & 80.93 & 52.37 \\
M3 & 80.44 & 55.23 \\
C1 & 76.14 & 51.26 \\
C2 (30-bus) & 76.95 & 52.29 \\
C2 (118-bus) & 78.72 & 53.80 \\
\hline
\end{tabular}
\end{table}

\subsection{Modified IEEE 30-Bus System}

Based on the prediction data mentioned above, the integrated framework in Fig.~\ref{fig:framework} and the proposed data-driven RO method in Algorithm~\ref{alg:RUC} are used for the robust UC problem. In the modified IEEE 30-bus system, there are six controllable generators and four wind farms, where $\underline{P}_g = (100, 40, 0, 0, 0, 0)$ MW, $\overline{P}_g = (360, 140, 100, 100, 100, 100)$ MW, and $T_g^+ = T_g^- = 6$ h. Other parameters are in \cite{xie2024github}. 

\subsubsection{Benchmark}

In the benchmark case, the probability tolerance parameters are set as $\varepsilon = \delta =$ 5\%. We divide the training dataset into two subsets with $N_1$ and $N_2$ data, respectively. According to Theorem~\ref{thm:construction} and Theorem~\ref{thm:reconstruction}, $N_2$ should be at least $59$. We set $N_1 = 212$ and $N_2 = 124$. 

Based on the performance on the validation dataset, we use the first three components of PCA in the surrogate model for weight optimization. The MLP neural network has two hidden layers, and each layer has 16 units. The neural network is trained using the Adam algorithm and the learning rate is 0.001. To help prevent overfitting, we adopt the L2 regularization technique. After about 1000 epochs, the neural network is trained. 

\begin{figure}[!t]
   \centering
   \includegraphics[width=0.48\textwidth, trim= 3mm 3mm 3mm 3mm, clip]{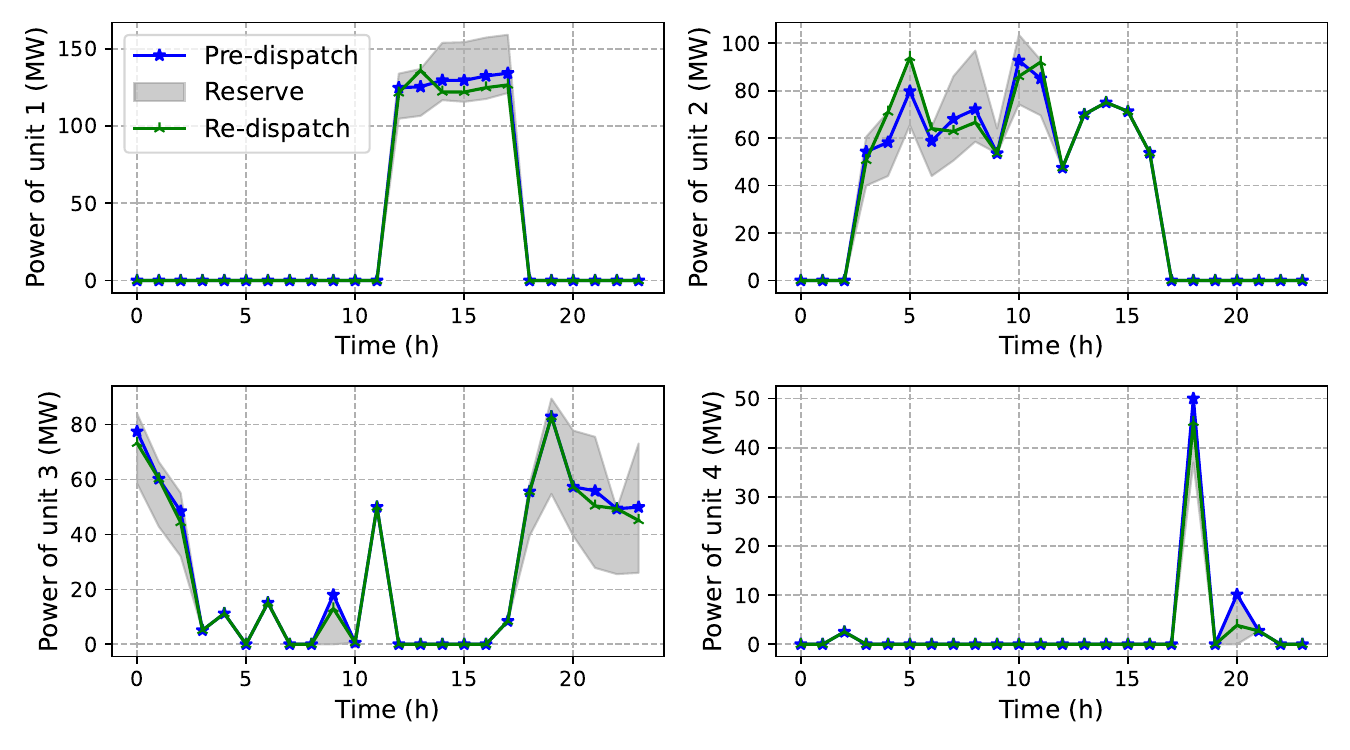}
   \caption{Power outputs of controllable generators.}
   \label{fig:unit}
\end{figure}

The test day for UC is from the test dataset. The MILP problem for weight optimization is solved in 0.02 s and the result is $w = (0.28, 0.23, 0.49)$. The weight for M3 is the highest, reflecting that the sub-profiles are valuable in power forecasting and robust UC, which is consistent with the results in TABLE~\ref{tab:load} and \cite{wang2018ensemble}. With the optimized weight, the UC results are obtained after 121 s, and the optimal value is \$89725. The power outputs of some controllable generators are depicted in Fig.~\ref{fig:unit}, including the pre-dispatch power, the reserve region, and the re-dispatch power in the worst-case scenario. The re-dispatch power is always within the reserve region.

\subsubsection{Comparison}

We use out-of-sample tests to examine the performance of the proposed method under uncertainty. Meanwhile, we compare the proposed method with the following alternative methods:
\begin{itemize}
    \item SP: The traditional chance-constrained SP method using the estimated distribution based on historical data.
    \item RO1: The traditional data-driven RO method with an ellipsoidal uncertainty set that includes all data.
    \item RO2: The same method as RO1 except that the ellipsoidal uncertainty set includes $1 - \varepsilon$ proportion of the data.
    \item P1: The data-driven RO method using the ellipsoidal uncertainty set with statistical guarantees.
    \item P2: The data-driven RO method using Algorithm~\ref{alg:RUC}, where the weight $w$ is optimized to minimize the MSE.
    \item PSO: Similar to the proposed method, but the weight is optimized using the particle swarm optimization (PSO) algorithm instead of the surrogate model. The number of particles is $10$, the inertia weight is $0.5$, the cognitive constant is $1.0$, the social constant is $1.5$, and the computation time limit is $10000$ s.
    \item EXACT: An adapted version of the MILP-based exact method for the predict-and-optimize approach in \cite{dias2025application}, whose details will be introduced later.
\end{itemize}
For clarity, the settings of these methods are compared in TABLE~\ref{tab:comparison}, with results listed in TABLE~\ref{tab:UC-30}. The feasible rate is calculated using the 100 forecast error samples on the test dataset. 
The total cost is calculated using the actual observed data of uncertainty on the test days.

\begin{table}[!t]
\renewcommand{\arraystretch}{1.3}
\caption{Settings of Unit Commitment Methods for Comparison}
\label{tab:comparison}
\centering
\begin{tabular}{cccc}
\hline
Method & Statistical & Integrated forecasting & Uncertainty set \\
& guarantee & and optimization & reconstruction \\
\hline
SP & $\times$ & $\times$ & $\times$ \\
RO1 & $\times$ & $\times$ & $\times$ \\
RO2 & $\times$ & $\times$ & $\times$ \\
P1 & \checkmark & \checkmark & $\times$ \\
P2 & \checkmark & $\times$ & \checkmark \\
PSO & \checkmark & \checkmark & \checkmark\\
EXACT & $\times$ & \checkmark & $\times$ \\
Proposed & \checkmark & \checkmark & \checkmark\\
\hline
\end{tabular}
\end{table}

\begin{table}[!t]
\renewcommand{\arraystretch}{1.3}
\caption{Unit Commitment Results of Different Methods in the Modified IEEE 30-Bus System}
\label{tab:UC-30}
\centering
\begin{tabular}{ccccc}
\hline
Method & Objective (\$) & Feasible rate & Test total cost (\$) & Time (s) \\
\hline
SP & 84832 & 88\% & 82985 & 218 \\
RO1 & 106810 & 100\% & 92652 & 143 \\
RO2 & 97350 & 97\% & 90468 & 94 \\
P1 & 97848 & 98\% & 89149 & 124 \\
P2 & 90122 & 98\% & 88318 & 147 \\
PSO & 89111 & 98\% & 88154 & 10000 \\
Proposed & 89725 & 98\% & 88243 & 121 \\
\hline
\end{tabular}
\end{table}

As TABLE~\ref{tab:UC-30} shows, SP has the lowest objective value and test total cost. However, its test feasible rate is $88$\%, much lower than the desired threshold $95$\%, which shows that SP lacks robustness. The other five methods are RO-based and their test feasible rates all meet the requirement. However, RO1 is rather conservative and has the highest objective value and test total cost. RO2 is a traditional data-driven RO method with no statistical guarantees, which shows its theoretical limitation. P1, P2, and the proposed method have statistical guarantees, and the proposed method achieves the lowest optimal objective and test total cost among them. Using uncertainty set reconstruction, the objective decreases $8.30$\%. The integrated forecasting and optimization framework also contributes to improving the objective.

To validate the effectiveness of the surrogate model, we compare the proposed method to a method that directly optimizes the weight. Since the impact of the weight on the evaluated cost of the UC strategy cannot be expressed as a compact function, we choose the PSO algorithm, which is derivative-free and highly flexible \cite{xie2018operationally}. As shown in TABLE~\ref{tab:UC-30}, the PSO algorithm uses a much longer time ($10000$ s) compared to the proposed method ($121$ s), achieving a reduction in the objective value by $0.7\%$ and a decrease in the test total cost by $0.1\%$. Thus, the proposed method substantially reduces the computation time to an acceptable duration for UC, with only a mild sacrifice in performance. This underscores the effectiveness of the proposed surrogate model in balancing computational efficiency and performance.

The EXACT method \cite{dias2025application} solves problem \eqref{eq:closed-loop}:
\begin{subequations}
\label{eq:closed-loop}
\begin{align}
    \label{eq:closed-loop-1}
    & \min_{w, \hat{\mathbf{U}}(w),x^*(w)}~ \frac{1}{N} \sum_{n = 1}^N \left( f(x_n^*(w)) + \min_{y \in \mathcal{Y}(x_n^*(w), \mathbf{U}_n)} h(y) \right) \\
    \label{eq:closed-loop-2}
    & \mbox{s.t.}~ \hat{\mathbf{U}}_n(w) = \sum\nolimits_{m \in \mathcal{M}} w^{(m)} \hat{\mathbf{U}}_n^{(m)}, \forall n, \\
    \label{eq:closed-loop-3}
    & x_n^*(w) \in \underset{x \in \mathcal{X}_n}{\operatorname{arg\,min}} \left\{ f(x) + \min_{y \in \mathcal{Y}(x, \hat{\mathbf{U}}_n(w))} h(y) \right\}, \forall n.
\end{align}
\end{subequations}
The predictions are combined using the weight $w = (w^{(m)}; m \in \mathcal{M})$ in \eqref{eq:closed-loop-2}. Based on the combined predictions, UC strategies are optimized in \eqref{eq:closed-loop-3}. The objective \eqref{eq:closed-loop-1} minimizes the average cost of UC strategies under observed uncertainty realizations.

The EXACT method uses the KKT conditions of the optimization problem in \eqref{eq:closed-loop-3} to equivalently reformulate problem \eqref{eq:closed-loop} into a single-level MILP. However, for the UC problem considered in this paper, there are binary variables in $x$, i.e., the indicator variables for the on/startup/shutdown states of generators, so the KKT conditions do not imply optimality. If we take a step back and relax the binary variables to continuous values in $[0, 1]$, the equivalent MILP for problem \eqref{eq:closed-loop} can be formulated. However, in the benchmark case with $N = 336$ training samples, the equivalent MILP deduced by KKT conditions contains over one million binary variables, and the Gurobi solver runs out of memory before finding a feasible solution. Even in a reduced dataset with $N = 40$ training data, the solver fails to converge within $10000$ s. Therefore, the EXACT method is not suitable for UC problems due to its computational complexity and limitations in handling binary variables in the original optimization problem.

Therefore, the above method comparison verifies that the proposed method is effective. It outperforms other methods considering the comprehensive performance of robustness, optimality, and computational efficiency. Moreover, the proposed weight optimization and uncertainty set reconstruction processes help improve the performance.

\begin{figure}[!t]
   \centering
   \includegraphics[width=0.36\textwidth, trim= 5mm 15mm 15mm 25mm, clip]{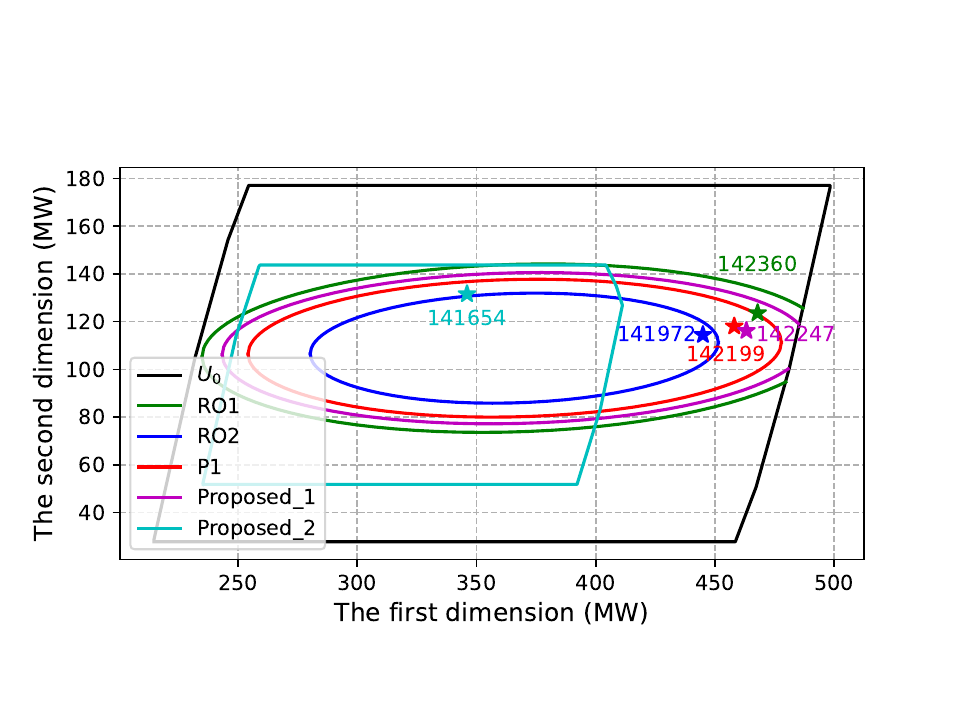}
   \caption{Projections of uncertainty sets onto a two-dimensional plane in the case of two random loads, where the projections of worst-case scenarios within these uncertainty sets are marked by stars, with the costs (\$) labeled.}
   \label{fig:projection}
\end{figure}

To visualize different uncertainty sets, we examine a specific case involving two random loads. 
The uncertainty sets are projected onto a two-dimensional plane. The projections of the uncertainty sets are drawn in Fig.~\ref{fig:projection}. To emphasize the impact of uncertainty set construction and reconstruction, the weight $w$ is fixed to that of C1. 
The black polygon is the projection of $\mathcal{U}_0$. 
All other uncertainty sets are the intersections of $\mathcal{U}_0$ and ellipsoids or polyhedrons. The projected uncertainty set of RO1 is framed by an ellipse that includes all the data points. RO2's ellipse has the same shape and center but only contains $95$\% data points. RO2 does not have a statistical guarantee for its out-of-sample performance. P1 has a larger uncertainty set than RO2 to maintain the statistical guarantees, but P1 does not reconstruct the uncertainty set to decrease the conservativeness. Proposed\_1 and Proposed\_2 are the projections of the first and second uncertainty sets of the proposed method ($\mathcal{U}_1'$ and $\mathcal{U}_2$ in Algorithm~\ref{alg:RUC}), respectively. Proposed\_1 is framed by an ellipse similar to that of P1. 
After the reconstruction, Proposed\_2 excludes the right part of the ellipse, thereby omitting some high-cost scenarios in the UC problem. At the same time, Proposed\_2 includes additional regions with relatively low costs, ensuring that the data points within it are enough for the statistical guarantees. This approach effectively reduces conservativeness while maintaining statistical guarantees.

\subsubsection{Sensitivity Analysis}

\begin{figure}[!t]
   \centering
   \includegraphics[width=0.48\textwidth, trim= 10mm 5mm 6mm 3mm, clip]{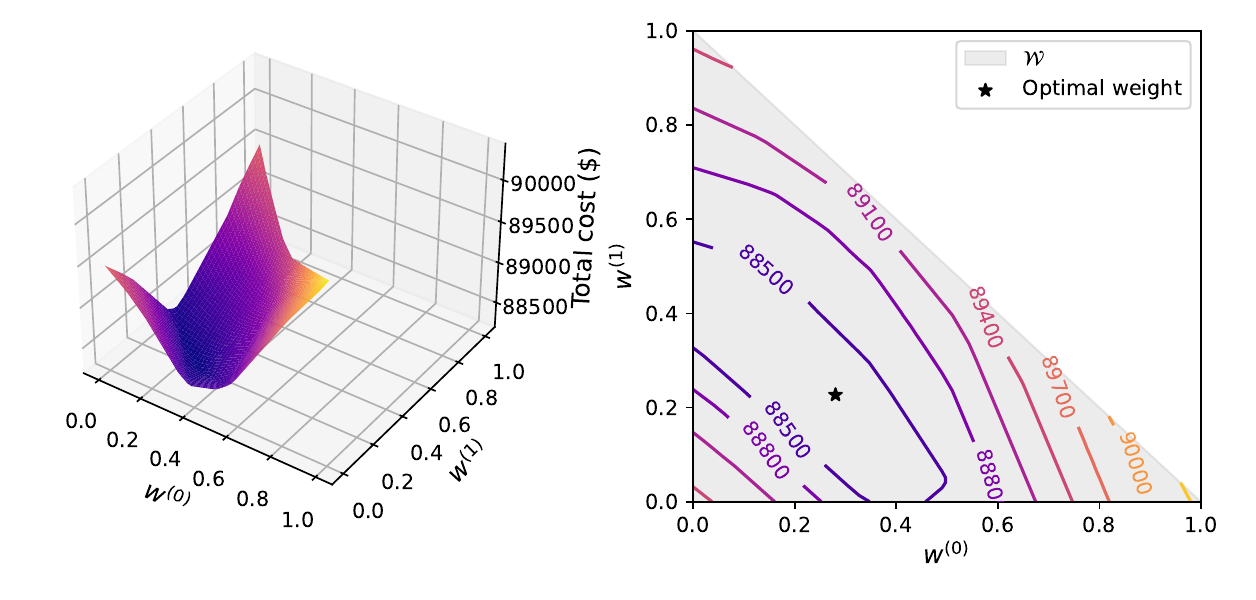}
   \caption{Total costs under different weights. (Because $w^{(2)} = 1 - w^{(0)} - w^{(1)}$, we only draw the relationship between the total cost, $w^{(0)}$, and $w^{(1)}$, where the weight satisfies $w^{(0)} \geq 0$, $w^{(1)} \geq 0$, and $w^{(0)} + w^{(1)} \leq 1$.)}
   \label{fig:weight}
\end{figure}

We investigate the impacts of weight $w$, the number of data points $N_2$, and the probability tolerance parameters $\varepsilon$ and $\delta$ for sensitivity analysis. The total costs under different values of weight $w$ are depicted in Fig.~\ref{fig:weight}. As Fig.~\ref{fig:weight} shows, there is a valley in the middle of the weight's feasible region, which means that the combination of the three kinds of predictions has the best performance. The optimal weight computed by the MILP problem of the surrogate model lies in the center region of the valley and achieves a low total cost, showing the effectiveness of the surrogate model.

\begin{figure}[!t]
   \centering
   \includegraphics[width=0.48\textwidth, trim= 3mm 5mm 0mm 3mm, clip]{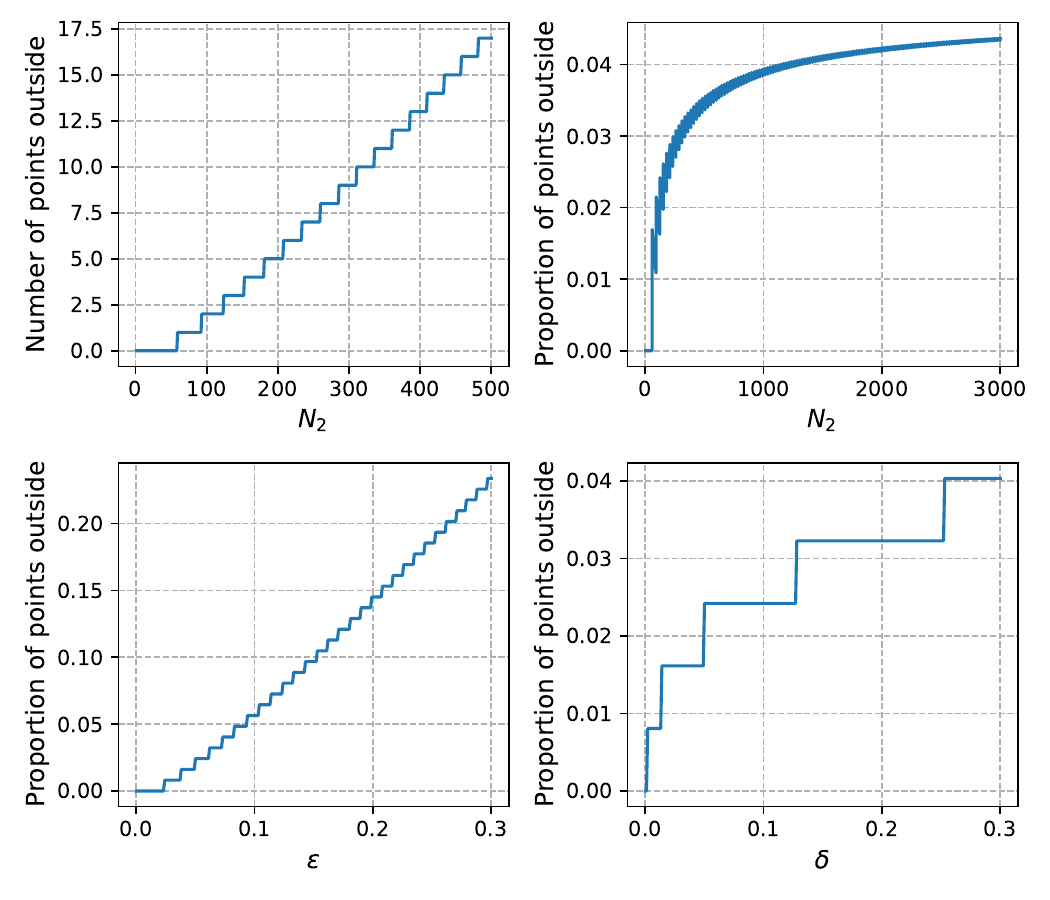}
   \caption{Proportion of points outside the proposed uncertainty set under different settings.}
   \label{fig:Nstar}
\end{figure}

To illustrate the relationship between statistical guarantees and the parameters $\varepsilon$ and $\delta$, we examine the number and proportion of the $N_2$ data points that fall outside the uncertainty set as $N_2$, $\varepsilon$, and $\delta$ vary.
In the benchmark case, $N_2 = 124$ and $\varepsilon = \delta = 5$\%. 
In Fig.~\ref{fig:Nstar}, as $N_2$ increases, the number of points outside steps up because it must be an integer. The proportion of points outside becomes closer to $\varepsilon = 5$\% under a larger $N_2$, indicating that statistical guarantees are more easily maintained with a larger dataset. When $\varepsilon$ is increased while keeping $N_2$ and $\delta$ fixed, the proportion of points outside the uncertainty set rises but never exceeds $\varepsilon$. This suggests that more data points should be included in the uncertainty set than the specified proportion to achieve a statistical guarantee of out-of-sample performance. 
As $\delta$ increases, the proportion of points outside increases.

\begin{figure}[!t]
   \centering
   \includegraphics[width=0.48\textwidth, trim= 3mm 5mm 0mm 3mm, clip]{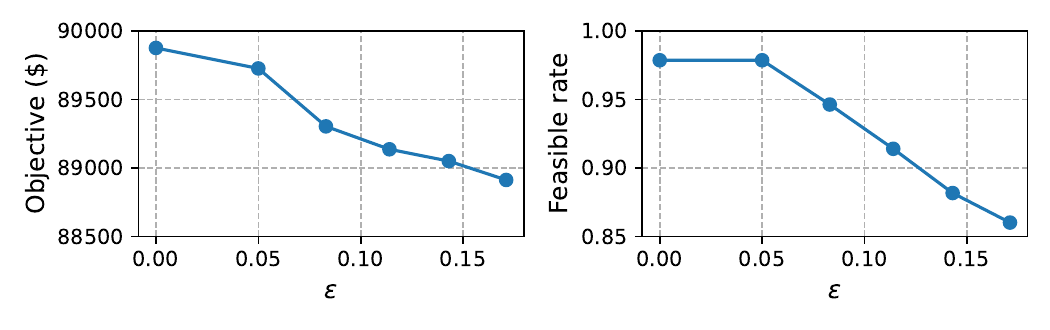}
   \caption{Results under different $\varepsilon$.}
   \label{fig:epsilon}
\end{figure}

Fig.~\ref{fig:epsilon} shows the results of the optimal value and the test feasible rate under different values of $\varepsilon$. As $\varepsilon$ increases, the robustness decreases, leading to a decline in both the objective value and the test feasible rate. When $\varepsilon \in [0.05, 0.17]$, the test feasible rate is larger than $1 - \varepsilon$, demonstrating the effectiveness of the probability guarantee. The test feasible rate decreases approximately linearly as $\varepsilon$ increases, when $\varepsilon \geq 0.05$. 

\begin{figure}[!t]
   \centering
   \includegraphics[width=0.48\textwidth, trim= 3mm 5mm 0mm 3mm, clip]{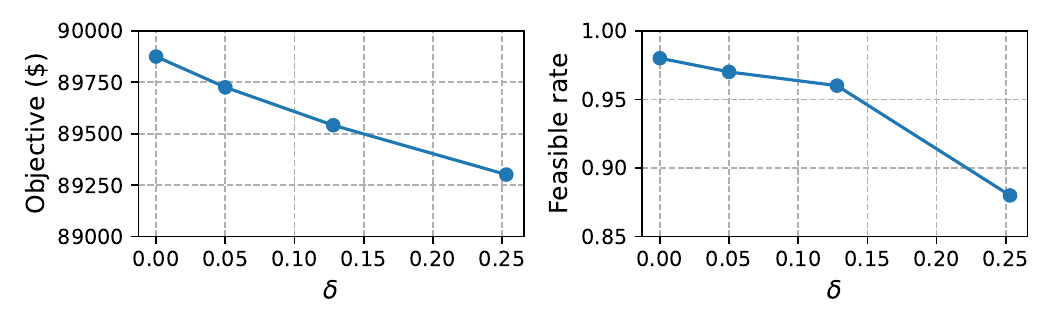}
   \caption{Results under different $\delta$.}
   \label{fig:delta}
\end{figure}

The effects of the probability threshold $\delta$ are illustrated in Fig.~\ref{fig:delta}. As $\delta$ increases, the confidence in the chance constraint $\Pr[u \in \mathcal{U}] \geq 1 - \varepsilon$ decreases, resulting in less conservative outcomes and decreases in both the objective value and the test feasible rate. When $\delta \leq 0.13$, the test feasible rate remains above $1 - \varepsilon = 95$\%. However, if $\delta$ becomes too large, the desired probability $1 - \varepsilon$ is no longer met in the test dataset.

Therefore, the parameters $\varepsilon$ and $\delta$ should be chosen carefully based on the decision makers' risk preference. Smaller values of $\varepsilon$ and $\delta$ result in more conservative strategies. We recommend using a small value for $\delta$ to maintain good out-of-sample performance, such as $\varepsilon = \delta = 5\%$.

\subsection{Modified IEEE 118-Bus System}

We use a modified IEEE 118-bus system with 54 controllable generators and four wind farms to show the scalability. The parameter settings are $\varepsilon = \delta =$ 5\%, $N_1 = 212$, and $N_2 = 124$. More details can be found in \cite{xie2024github}. The UC results under different methods are shown in TABLE~\ref{tab:UC-118}. The method relationships are similar to those in the modified IEEE 30-bus system. In addition, the computation time increases but is still acceptable for day-ahead UC. We also test the computational efficiency under different numbers of random loads in the modified IEEE 118-bus system, as shown in Table~\ref{tab:time}. These results verify the scalability of the proposed method.

\begin{table}[!t]
\renewcommand{\arraystretch}{1.3}
\caption{Unit Commitment Results of Different Methods in the Modified IEEE 118-Bus System}
\label{tab:UC-118}
\centering
\begin{tabular}{ccccc}
\hline
Method & Objective (\$) & Feasible rate & Test total cost (\$) & Time (s) \\
\hline
SP & 2061915 & 84\% & 2055123 & 1638 \\
RO1 & 2096428 & 100\% & 2069461 & 256 \\
RO2 & 2086813 & 100\% & 2064276 & 634 \\
P1 & 2092717 & 100\% & 2068213 & 1304 \\
P2 & 2066737 & 100\% & 2058488 & 1076 \\
Proposed & 2065591 & 100\% & 2056996 & 763 \\
\hline
\end{tabular}
\end{table}

\begin{table}[!t]
\renewcommand{\arraystretch}{1.3}
\caption{Computational Efficiency in the Modified IEEE 118-Bus System}
\label{tab:time}
\centering
\begin{tabular}{ccccc}
\hline
Number of random loads & 25 & 21 & 17 & 13 \\
\hline
Number of iterations & 23 & 21 & 21 & 18 \\
Computation time (s) & 763 & 600 & 551 & 391 \\
\hline
\end{tabular}
\end{table}

\section{Conclusion}
\label{sec:conclusion}

To enhance out-of-sample performance and ensure robustness, this paper develops a new data-driven two-stage robust UC method. The proposed integrated forecasting and optimization framework combines different predictions using weights optimized based on performance outcomes. A surrogate model is established to accelerate the weight optimization process. In the two-stage robust UC, the uncertainty set is constructed from data to have statistical guarantees, and it is then reconstructed using the information from the optimization problem to reduce conservativeness. Comparative analysis on modified IEEE 30-bus and 118-bus systems demonstrates that the proposed method surpasses traditional SP and RO methods in balancing robustness with out-of-sample performance. The case studies also show that the computational complexity of the proposed method is comparable to that of traditional RO methods, making it scalable. Future work could explore multi-stage robust UC with statistical guarantees.

\appendices

\section{Proof of Lemma~\ref{lemma:RUC}}
\label{appendix:lemma-RUC}
\setcounter{equation}{0}  
\renewcommand{\theequation}{A.\arabic{equation}}

1) Prove that $O = O_{x^*} \leq O_{x_\mathcal{U}^*} \leq O_\mathcal{U}$ when $\Pr[u \in \mathcal{U}] \geq 1 - \varepsilon$:

Because $(x^*, \eta^*)$ is optimal in \eqref{eq:UC-original}, we have $O = O_{x^*}$ and 
\begin{align}
    O & = f(x^*) + \min \left\{ \eta ~\middle|~ \Pr \left[ \min_{y \in \mathcal{Y}(x^*, u)} h(y) \leq \eta \right] \geq 1 - \varepsilon \right\} \nonumber \\
    & \leq f(x_\mathcal{U}^*) + \min \left\{ \eta ~\middle|~ \Pr \left[ \min_{y \in \mathcal{Y}(x_\mathcal{U}^*, u)} h(y) \leq \eta \right] \geq 1 - \varepsilon \right\} \nonumber \\
    & = O_{x_\mathcal{U}^*}. \nonumber
\end{align}
Let 
\begin{align}
    \eta_\mathcal{U}^* = \max\nolimits_{u \in \mathcal{U}} \min\nolimits_{y \in \mathcal{Y}(x_\mathcal{U}^*, u)} h(y). \nonumber
\end{align}
Then 
\begin{align}
    \Pr \left[ \min\nolimits_{y \in \mathcal{Y}(x_\mathcal{U}^*, u)} h(y) \leq \eta_\mathcal{U}^* \right] \geq \Pr[u \in \mathcal{U}] \geq 1 - \varepsilon. \nonumber
\end{align}
Therefore,
\begin{align}
    O_{x_\mathcal{U}^*} & = f(x_\mathcal{U}^*) + \min \left\{ \eta ~\middle|~ \Pr \left[ \min_{y \in \mathcal{Y}(x_\mathcal{U}^*, u)} h(y) \leq \eta \right] \geq 1 - \varepsilon \right\} \nonumber \\
    & \leq f(x_\mathcal{U}^*) + \eta_\mathcal{U}^* = f(x_\mathcal{U}^*) + \max_{u \in \mathcal{U}} \min_{y \in \mathcal{Y}(x_\mathcal{U}^*, u)} h(y) = O_\mathcal{U}, \nonumber
\end{align}
where the last equation follows from the optimality of $x_\mathcal{U}^*$ in \eqref{eq:RUC-original}. Hence, $O = O_{x^*} \leq O_{x_\mathcal{U}^*} \leq O_\mathcal{U}$.

2) Prove that $\Pr[u \in \mathcal{U}^*] \geq 1 - \varepsilon$ and $O = O_{\mathcal{U}^*}$:

By the feasibility of $(x^*, \eta^*)$ in \eqref{eq:UC-original}, we have
\begin{align}
    \Pr[u \in \mathcal{U}^*] = \Pr \left[ \min\nolimits_{y \in \mathcal{Y}(x^*, u)} h(y) \leq \eta^* \right] \geq 1 - \varepsilon, \nonumber
\end{align}
and then $O \leq O_{\mathcal{U}^*}$ by the previous conclusion. According to the definition of $\mathcal{U}^*$ in \eqref{eq:U*}, $\min_{y \in \mathcal{Y}(x^*, u)} h(y) \leq \eta^*$ for any $u \in \mathcal{U}^*$. Therefore,
\begin{align}
    O_{\mathcal{U}^*} & = \min\nolimits_{x \in \mathcal{X}} \left\{ f(x) + \max\nolimits_{u \in \mathcal{U}^*} \min\nolimits_{y \in \mathcal{Y}(x, u)} h(y) \right\} \nonumber \\
    & \leq f(x^*) + \max\nolimits_{u \in \mathcal{U}^*} \min\nolimits_{y \in \mathcal{Y}(x^*, u)} h(y) \nonumber \\
    & \leq f(x^*) + \eta^* = O. \nonumber
\end{align}
This completes the proof.

\section{Proof of Lemma~\ref{lemma:ramp}}
\label{appendix:lemma-ramp}
\setcounter{equation}{0}  
\renewcommand{\theequation}{B.\arabic{equation}}

According to the definitions of $\mathcal{X}$ and $\mathcal{Y}(x, u)$, we have \eqref{eq:X-l}, \eqref{eq:X-m}, and \eqref{eq:Y-c}. Constraint \eqref{eq:Y-c} implies that
\begin{align}
    & p_{g t} + p_{g t}^+ - p_{g t}^- \leq p_{g t} + p_{g t}^+ \leq p_{g t} + r_{g t}^+, \forall g \in \mathcal{G}, \forall t \in \mathcal{T}, \nonumber \\
    & p_{g t} + p_{g t}^+ - p_{g t}^- \geq p_{g t} - p_{g t}^- \geq p_{g t} - r_{g t}^-, \forall g \in \mathcal{G}, \forall t \in \mathcal{T}. \nonumber
\end{align}
Thus, by \eqref{eq:X-l}, 
\begin{align}
    & (p_{g t}+ p_{g t}^+ - p_{g t}^-) - (p_{g (t-1)}+ p_{g (t-1)}^+ - p_{g (t-1)}^-) \nonumber \\
    \leq\, & (p_{g t} + r_{g t}^+) - (p_{g (t-1)} - r_{g (t-1)}^-) \nonumber \\
    \leq\, & K_g^+ \theta_{g (t-1)} + K_g^U \theta_{g t}^+, \forall g \in \mathcal{G}, 2 \leq t \leq T, \nonumber
\end{align}
which shows \eqref{eq:ramp-re-upper}. Constraint \eqref{eq:ramp-re-lower} follows similarly.

\section{Proof of Theorem~\ref{thm:construction}}
\label{appendix:thm-construction}
\setcounter{equation}{0}  
\renewcommand{\theequation}{C.\arabic{equation}}

The function $a(e)$ is a continuous random variable, where $e$ follows the distribution $\mathbb{P}$. Therefore, the assumptions of Theorem 1 in \cite{hong2021learning} are satisfied, and then its conclusion holds, i.e., $\mathcal{E}$ satisfies the statistical guarantee \eqref{eq:guarantee-error}. 

Because $u \in \mathcal{U}_0$ must hold, we have $u - \hat{u} = e \in \mathcal{E} \iff u \in \mathcal{U}_1$. Combine it with Lemma~\ref{lemma:RUC}, then
\begin{align}
    \Pr[e \in \mathcal{E}] \geq 1 - \varepsilon & \iff \Pr[u \in \mathcal{U}_1] \geq 1 - \varepsilon \nonumber \\
    & \implies O \leq O_{x_0} \leq O_{\mathcal{U}_1}. \nonumber
\end{align}
Hence, \eqref{eq:guarantee-construction} follows from \eqref{eq:guarantee-error}, which completes the proof.

\section{Proof of Theorem~\ref{thm:reconstruction}}
\label{appendix:thm-reconstruction}
\setcounter{equation}{0}  
\renewcommand{\theequation}{D.\arabic{equation}}

1) Prove \eqref{eq:guarantee-reconstruction-1}:

According to linear programming theory \cite{boyd2004convex}, the function $b(\hat{u} + e)$ in $e$ is continuous on the closed set 
\begin{align}
    & \left\{ e ~\middle|~ \min\nolimits_{y \in \mathcal{Y}(x_0, \hat{u} + e)} h(y) < + \infty \right\} \nonumber \\
    = & \left\{e ~\middle|~ \mathcal{Y}(x_0, \hat{u} + e) \neq \varnothing\right\} \nonumber \\
    = & \left\{e ~\middle|~ \exists y, ~\mbox{s.t.}~ A y \geq B x_0 + D(\hat{u} + e) + E\right\}. \nonumber
\end{align}
When $\beta = +\infty$, $\mathcal{U}_2 = \mathcal{U}_0$ and \eqref{eq:guarantee-reconstruction-1} holds. For the finite case, Lemma 3 and Theorem 1 in \cite{hong2021learning} can be applied to $b(\hat{u} + e)$ to obtain \eqref{eq:guarantee-reconstruction-1}.

2) Prove \eqref{eq:guarantee-reconstruction-2}: 

Similar to the proof of Theorem~\ref{thm:construction}, the statistical guarantee \eqref{eq:guarantee-reconstruction-1} implies
\begin{align}
    \mathbb{P}^N[O \leq O_{x_1} \leq O_{\mathcal{U}_2}] \geq 1 - \delta. \nonumber
\end{align}
In addition,
\begin{align}
    O_{\mathcal{U}_2} & = \min\nolimits_{x \in \mathcal{X}} \left\{ f(x) + \max\nolimits_{u \in \mathcal{U}_2} \min\nolimits_{y \in \mathcal{Y}(x, u)} h(y) \right\} \nonumber \\
    & \leq f(x_0) + \max\nolimits_{u \in \mathcal{U}_2} \min\nolimits_{y \in \mathcal{Y}(x_0, u)} h(y) \leq f(x_0) + \beta, \nonumber
\end{align}
where the last inequality follows from \eqref{eq:U2}, so \eqref{eq:guarantee-reconstruction-2} holds.

\bibliographystyle{IEEEtran}
\bibliography{mybib}

\begin{thebibliography}{10}
\providecommand{\url}[1]{#1}
\csname url@samestyle\endcsname
\providecommand{\newblock}{\relax}
\providecommand{\bibinfo}[2]{#2}
\providecommand{\BIBentrySTDinterwordspacing}{\spaceskip=0pt\relax}
\providecommand{\BIBentryALTinterwordstretchfactor}{4}
\providecommand{\BIBentryALTinterwordspacing}{\spaceskip=\fontdimen2\font plus
\BIBentryALTinterwordstretchfactor\fontdimen3\font minus \fontdimen4\font\relax}
\providecommand{\BIBforeignlanguage}[2]{{%
\expandafter\ifx\csname l@#1\endcsname\relax
\typeout{** WARNING: IEEEtran.bst: No hyphenation pattern has been}%
\typeout{** loaded for the language `#1'. Using the pattern for}%
\typeout{** the default language instead.}%
\else
\language=\csname l@#1\endcsname
\fi
#2}}
\providecommand{\BIBdecl}{\relax}
\BIBdecl

\bibitem{paturet2020stochastic}
M.~Paturet, U.~Markovic, S.~Delikaraoglou, E.~Vrettos, P.~Aristidou, and G.~Hug, ``Stochastic unit commitment in low-inertia grids,'' \emph{IEEE Transactions on Power Systems}, vol.~35, no.~5, pp. 3448--3458, 2020.

\bibitem{xiong2022multi}
H.~Xiong, Y.~Shi, Z.~Chen, C.~Guo, and Y.~Ding, ``Multi-stage robust dynamic unit commitment based on pre-extended-fast robust dual dynamic programming,'' \emph{IEEE Transactions on Power Systems}, vol.~38, no.~3, pp. 2411--2422, 2022.

\bibitem{zhou2023partial}
B.~Zhou, J.~Fang, X.~Ai, Y.~Zhang, W.~Yao, Z.~Chen, and J.~Wen, ``Partial-dimensional correlation-aided convex-hull uncertainty set for robust unit commitment,'' \emph{IEEE Transactions on Power Systems}, vol.~38, no.~3, pp. 2434--2446, 2023.

\bibitem{ju2023two}
C.~Ju, T.~Ding, W.~Jia, C.~Mu, H.~Zhang, and Y.~Sun, ``Two-stage robust unit commitment with the cascade hydropower stations retrofitted with pump stations,'' \emph{Applied Energy}, vol. 334, p. 120675, 2023.

\bibitem{wang2024two}
W.~Wang, A.~Danandeh, B.~Buckley, and B.~Zeng, ``Two-stage robust unit commitment problem with complex temperature and demand uncertainties,'' \emph{IEEE Transactions on Power Systems}, vol.~39, no.~1, pp. 909--920, 2024.

\bibitem{haghighat2024robust}
H.~Haghighat, W.~Wang, and B.~Zeng, ``Robust unit commitment with decision-dependent uncertain demand and time-of-use pricing,'' \emph{IEEE Transactions on Power Systems}, vol.~39, no.~2, pp. 2854--2865, 2024.

\bibitem{wang2022distributionally}
S.~Wang, C.~Zhao, L.~Fan, and R.~Bo, ``Distributionally robust unit commitment with flexible generation resources considering renewable energy uncertainty,'' \emph{IEEE Transactions on Power Systems}, vol.~37, no.~6, pp. 4179--4190, 2022.

\bibitem{zheng2023day}
X.~Zheng, B.~Zhou, X.~Wang, B.~Zeng, J.~Zhu, H.~Chen, and W.~Zheng, ``Day-ahead network-constrained unit commitment considering distributional robustness and intraday discreteness: A sparse solution approach,'' \emph{Journal of Modern Power Systems and Clean Energy}, vol.~11, no.~2, pp. 489--501, 2023.

\bibitem{zhou2023distributionally}
A.~Zhou, M.~Yang, X.~Zheng, and S.~Yin, ``Distributionally robust unit commitment considering unimodality-skewness information of wind power uncertainty,'' \emph{IEEE Transactions on Power Systems}, vol.~38, no.~6, pp. 5420--5431, 2023.

\bibitem{liu2024modeling}
L.~Liu, Z.~Hu, Y.~Wen, and Y.~Ma, ``Modeling of frequency security constraints and quantification of frequency control reserve capacities for unit commitment,'' \emph{IEEE Transactions on Power Systems}, vol.~39, no.~1, pp. 2080--2092, 2024.

\bibitem{gao2023finite}
R.~Gao, ``Finite-sample guarantees for wasserstein distributionally robust optimization: Breaking the curse of dimensionality,'' \emph{Operations Research}, vol.~71, no.~6, pp. 2291--2306, 2023.

\bibitem{xie2022sizing}
R.~Xie, W.~Wei, M.~Shahidehpour, Q.~Wu, and S.~Mei, ``Sizing renewable generation and energy storage in stand-alone microgrids considering distributionally robust shortfall risk,'' \emph{IEEE Transactions on Power Systems}, vol.~37, no.~5, pp. 4054--4066, 2022.

\bibitem{hong2021learning}
L.~J. Hong, Z.~Huang, and H.~Lam, ``Learning-based robust optimization: Procedures and statistical guarantees,'' \emph{Management Science}, vol.~67, no.~6, pp. 3447--3467, 2021.

\bibitem{lu2024sample}
C.~Lu, N.~Gu, W.~Jiang, and C.~Wu, ``Sample-adaptive robust economic dispatch with statistically feasible guarantees,'' \emph{IEEE Transactions on Power Systems}, vol.~39, no.~1, pp. 779--793, 2024.

\bibitem{jiang2023robust}
W.~Jiang, C.~Lu, and C.~Wu, ``Robust scheduling of thermostatically controlled loads with statistically feasible guarantees,'' \emph{IEEE Transactions on Smart Grid}, vol.~14, no.~5, pp. 3561--3572, 2023.

\bibitem{liang2024joint}
J.~Liang, W.~Jiang, C.~Lu, and C.~Wu, ``Joint chance-constrained unit commitment: Statistically feasible robust optimization with learning-to-optimize acceleration,'' \emph{IEEE Transactions on Power Systems}, 2024.

\bibitem{wang2018ensemble}
Y.~Wang, Q.~Chen, M.~Sun, C.~Kang, and Q.~Xia, ``An ensemble forecasting method for the aggregated load with subprofiles,'' \emph{IEEE Transactions on Smart Grid}, vol.~9, no.~4, pp. 3906--3908, 2018.

\bibitem{cui2022ensemble}
W.~Cui, C.~Wan, and Y.~Song, ``Ensemble deep learning-based non-crossing quantile regression for nonparametric probabilistic forecasting of wind power generation,'' \emph{IEEE Transactions on Power Systems}, vol.~38, no.~4, pp. 3163--3178, 2022.

\bibitem{wang2023risk}
J.~Wang, Y.~Zhou, Y.~Zhang, F.~Lin, and J.~Wang, ``Risk-averse optimal combining forecasts for renewable energy trading under cvar assessment of forecast errors,'' \emph{IEEE Transactions on Power Systems}, vol.~39, no.~1, pp. 2296--2309, 2023.

\bibitem{li2023adaptive}
M.~Li, M.~Yang, Y.~Yu, M.~Shahidehpour, and F.~Wen, ``Adaptive weighted combination approach for wind power forecast based on deep deterministic policy gradient method,'' \emph{IEEE Transactions on Power Systems}, vol.~39, no.~2, pp. 3075--3087, 2024.

\bibitem{su2024towards}
H.-Y. Su and C.-C. Lai, ``Towards improved load forecasting in smart grids: A robust deep ensemble learning framework,'' \emph{IEEE Transactions on Smart Grid}, 2024.

\bibitem{elmachtoub2022smart}
A.~N. Elmachtoub and P.~Grigas, ``Smart “predict, then optimize”,'' \emph{Management Science}, vol.~68, no.~1, pp. 9--26, 2022.

\bibitem{dias2025application}
J.~Dias~Garcia, A.~Street, T.~Homem-de Mello, and F.~D. Mu{\~n}oz, ``Application-driven learning: A closed-loop prediction and optimization approach applied to dynamic reserves and demand forecasting,'' \emph{Operations Research}, vol.~73, no.~1, pp. 22--39, 2025.

\bibitem{chen2022feature}
X.~Chen, Y.~Yang, Y.~Liu, and L.~Wu, ``Feature-driven economic improvement for network-constrained unit commitment: A closed-loop predict-and-optimize framework,'' \emph{IEEE Transactions on Power Systems}, vol.~37, no.~4, pp. 3104--3118, 2022.

\bibitem{sang2022safety}
L.~Sang, Y.~Xu, H.~Long, and W.~Wu, ``Safety-aware semi-end-to-end coordinated decision model for voltage regulation in active distribution network,'' \emph{IEEE Transactions on Smart Grid}, vol.~14, no.~3, pp. 1814--1826, 2022.

\bibitem{wu2024novel}
H.~Wu, D.~Ke, L.~Song, S.~Liao, J.~Xu, Y.~Sun, and K.~Fang, ``A novel stochastic unit commitment characterized by closed-loop forecast-and-decision for wind integrated power systems,'' \emph{IEEE Transactions on Power Systems}, vol.~39, no.~2, pp. 2570--2586, 2024.

\bibitem{shen2024predict}
F.~Shen, Y.~Cao, M.~Shahidehpour, X.~Xu, C.~Wang, J.~Wang, and S.~Zhai, ``Predict-and-optimize model for day-ahead inertia prediction using distributionally robust unit commitment with renewable energy sources,'' \emph{IEEE Transactions on Power Systems}, 2024.

\bibitem{zhou2024load}
Y.~Zhou, Q.~Wen, J.~Song, X.~Cui, and Y.~Wang, ``Load data valuation in multi-energy systems: An end-to-end approach,'' \emph{IEEE Transactions on Smart Grid}, 2024.

\bibitem{zhuang2025weighted}
Y.~Zhuang, L.~Cheng, C.~Wan, R.~Xie, N.~Qi, and Y.~Chen, ``A weighted predict-and-optimize framework for power system operation considering varying impacts of uncertainty,'' \emph{arXiv preprint arXiv:2503.11001}, 2025.

\bibitem{abdi2010principal}
H.~Abdi and L.~J. Williams, ``Principal component analysis,'' \emph{Wiley interdisciplinary reviews: computational statistics}, vol.~2, no.~4, pp. 433--459, 2010.

\bibitem{anderson2020strong}
R.~Anderson, J.~Huchette, W.~Ma, C.~Tjandraatmadja, and J.~P. Vielma, ``Strong mixed-integer programming formulations for trained neural networks,'' \emph{Mathematical Programming}, vol. 183, pp. 3--39, 2020.

\bibitem{leshno1993multilayer}
M.~Leshno, V.~Y. Lin, A.~Pinkus, and S.~Schocken, ``Multilayer feedforward networks with a nonpolynomial activation function can approximate any function,'' \emph{Neural networks}, vol.~6, no.~6, pp. 861--867, 1993.

\bibitem{lu2017expressive}
Z.~Lu, H.~Pu, F.~Wang, Z.~Hu, and L.~Wang, ``The expressive power of neural networks: A view from the width,'' \emph{Advances in neural information processing systems}, vol.~30, 2017.

\bibitem{arroyo2000optimal}
J.~M. Arroyo and A.~J. Conejo, ``Optimal response of a thermal unit to an electricity spot market,'' \emph{IEEE Transactions on power systems}, vol.~15, no.~3, pp. 1098--1104, 2000.

\bibitem{moreira2024role}
A.~Moreira, B.~Fanzeres, P.~Silva, M.~Heleno, and A.~L.~M. Marcato, ``On the role of battery energy storage systems in the day-ahead contingency-constrained unit commitment problem under renewable penetration,'' \emph{Electric Power Systems Research}, vol. 235, p. 110856, 2024.

\bibitem{velloso2019two}
A.~Velloso, A.~Street, D.~Pozo, J.~M. Arroyo, and N.~G. Cobos, ``Two-stage robust unit commitment for co-optimized electricity markets: An adaptive data-driven approach for scenario-based uncertainty sets,'' \emph{IEEE Transactions on Sustainable Energy}, vol.~11, no.~2, pp. 958--969, 2019.

\bibitem{golestaneh2018ellipsoidal}
F.~Golestaneh, P.~Pinson, R.~Azizipanah-Abarghooee, and H.~B. Gooi, ``Ellipsoidal prediction regions for multivariate uncertainty characterization,'' \emph{IEEE Transactions on Power Systems}, vol.~33, no.~4, pp. 4519--4530, 2018.

\bibitem{zeng2013solving}
B.~Zeng and L.~Zhao, ``Solving two-stage robust optimization problems using a column-and-constraint generation method,'' \emph{Operations Research Letters}, vol.~41, no.~5, pp. 457--461, 2013.

\bibitem{grabner2023global}
M.~Grabner, Y.~Wang, Q.~Wen, B.~Bla{\v{z}}i{\v{c}}, and V.~{\v{S}}truc, ``A global modeling framework for load forecasting in distribution networks,'' \emph{IEEE Transactions on Smart Grid}, vol.~14, no.~6, pp. 4927--4941, 2023.

\bibitem{commission2012cer}
{Commission for Energy Regulation (CER)}, ``Cer smart metering project-electricity customer behaviour trial, 2009-2010,'' 2012.

\bibitem{staffell2016using}
I.~Staffell and S.~Pfenninger, ``Using bias-corrected reanalysis to simulate current and future wind power output,'' \emph{Energy}, vol. 114, pp. 1224--1239, 2016.

\bibitem{siami2019performance}
S.~Siami-Namini, N.~Tavakoli, and A.~S. Namin, ``The performance of lstm and bilstm in forecasting time series,'' in \emph{2019 IEEE International conference on big data (Big Data)}.\hskip 1em plus 0.5em minus 0.4em\relax IEEE, 2019, pp. 3285--3292.

\bibitem{fekri2022distributed}
M.~N. Fekri, K.~Grolinger, and S.~Mir, ``Distributed load forecasting using smart meter data: Federated learning with recurrent neural networks,'' \emph{International Journal of Electrical Power \& Energy Systems}, vol. 137, p. 107669, 2022.

\bibitem{xie2024github}
R.~Xie, ``data-driven-robust-unit-commitment,'' \url{https://github.com/xieruijx/data-driven-robust-unit-commitment}, 2025.

\bibitem{xie2018operationally}
R.~Xie, Y.~Chen, F.~Li, Z.~Wang, and S.~Mei, ``Operationally constrained optimal dispatch of multiple pulsed loads in an isolated microgrid,'' in \emph{2018 IEEE Power \& Energy Society General Meeting (PESGM)}.\hskip 1em plus 0.5em minus 0.4em\relax IEEE, 2018, pp. 1--5.

\bibitem{boyd2004convex}
S.~P. Boyd and L.~Vandenberghe, \emph{Convex optimization}.\hskip 1em plus 0.5em minus 0.4em\relax Cambridge university press, 2004.

\end{thebibliography}

\end{document}